\documentclass[3p,times,twocolumn]{elsarticle}

\usepackage{amsmath,amsfonts,amssymb}
\biboptions{authoryear}
\usepackage{mathrsfs}
\usepackage{ulem}
\usepackage{graphicx}
\usepackage{graphics,color}
\usepackage{float}
\usepackage{stfloats}
\usepackage{epsfig}
\usepackage{natbib}
\usepackage{epstopdf}
\usepackage{algorithm}
\usepackage{algpseudocode}
\usepackage{xcolor}
\usepackage{geometry} 
\usepackage[colorlinks=true,pdfstartview=FitH,linkcolor=blue,anchorcolor=violet,citecolor=magenta]{hyperref}

\newtheorem{lemma}{Lemma}
\newtheorem{theorem}{Theorem}
\newtheorem{corollary}{Corollary}
\newtheorem{definition}{Definition}
\newtheorem{remark}{Remark}
\newtheorem{example}{Example}
\newtheorem{assumption}{Assumption}
\newtheorem{notation}{Notation:}

\let\oldremark\remark
\renewcommand{\remark}{\oldremark\normalfont}

\let\oldtheorem\theorem
\renewcommand{\theorem}{\oldtheorem\normalfont}

\let\oldlemma\lemma
\renewcommand{\lemma}{\oldlemma\normalfont}

\let\olddefinition\definition
\renewcommand{\definition}{\olddefinition\normalfont}

\let\oldcorollary\corollary
\renewcommand{\corollary}{\oldcorollary\normalfont}

\let\oldexample\example
\renewcommand{\example}{\oldexample\normalfont}

\let\oldassumption\assumption
\renewcommand{\assumption}{\oldassumption\normalfont}

\let\oldnotation\notation
\renewcommand{\notation}{\oldnotation\normalfont}

\usepackage[labelsep=period]{caption}

\allowdisplaybreaks
\begin{document}
	\begin{frontmatter}
		\title{{\color{black}AW-EL-PINNs: A Multi-Task Learning Physics-Informed Neural Network for Euler-Lagrange Systems in Optimal Control Problems\tnoteref{t1}}}
		\tnotetext[t1]{This work was supported by the
			National Natural Science Foundation of China under Grant 62373310,
			and the Graduate Student
			Research Innovation Project of Southwest University under Grant SWUB25081.}

		\author[address1]{Chuandong Li}
		\ead{cdli@swu.edu.cn}
		
		\author[address1]{Runtian Zeng\corref{correspondingauthor}}
		\cortext[correspondingauthor]{Corresponding author}
		\ead{zzzrt312@email.swu.edu.cn}

		\address[address1]{College of Electronic and Information Engineering, Southwest University, Chongqing, 400715, China}
		
		\begin{abstract}
			This paper presents adaptive weighted Euler-Lagrange theorem combined physics-informed neural networks (AW-EL-PINNs) for solving Euler-Lagrange systems in optimal control problems. The framework systematically converts optimal control frameworks into two-point boundary value problems (TPBVPs) while establishing a multi-task learning paradigm through innovative integration of the Euler-Lagrange theorem with deep learning  architecture. An adaptive loss weighting mechanism  dynamically balances loss function components during training, decreasing tedious manual tuning of weighting the loss functions compared to the conventional physics-informed neural networks (PINNs). Based on  six numerical examples, it's clear  that AW-EL-PINNs achieve enhanced solution accuracy compared to baseline methods while maintaining   stability throughout the optimization process. These results highlight the framework’s  capability to improve precision and ensure  stability in solving   Euler-Lagrange systems in optimal control problems, offering  potential strategies for problems under  physical applications.
		\end{abstract}
		
		\begin{keyword}
			
			Euler-Lagrange theorem;
			Physics-informed neural networks; 
			Optimal control problems; 
			Two-point boundary value problems;
			Adaptive loss weighting

		\end{keyword}
		
	\end{frontmatter}

	\section{Introduction}\label{Sec:Intro.}
	The optimal control problems constitute a category of problems that focus on searching for the optimal strategies, which  maximize or minimize the performance indicator function \citep{kirk2004optimal}. Such  problems are commonly found in the modeling of fields such as  aerospace control \citep{schiassi2020physics,schiassi2022physics}, robot control \citep{ishihara2018optimal,dong2025barrier},  and automatic driving \citep{liu2023neurodynamic,li2024pontryagin}. Finding optimal strategies for optimal control problems is equivalent to solving the problems, and the Pontryagin minimum principle is a classic solution method. It transforms optimal control problems into two-point boundary value problems (TPBVPs) involving state variables and adjoint variables, thereby providing necessary conditions for optimality and enabling the determination of optimal solutions.
	
	However, obtaining analytical expressions for solutions is not a trivial task. Even for structurally simple problems, such as linear quadratic (LQ) optimal control problems, analytical solutions do not always exist.
	Therefore, the computational complexity of optimal control problems has driven significant interest in numerical solutions, particularly the gradient method and conjugate gradient method 
	\citep{lasdon2003conjugate}. These approaches fundamentally rely on the Euler-Lagrange theorem—a special case of the Pontryagin minimum principle applicable to unconstrained control problems.
	Specifically, when the control domain is unconstrained, the Hamiltonian function is smooth, and the optimal solution lies in the interior of the feasible domain, the  problems can reach a stationary point under the condition $\frac{\partial H}{\partial u}=0$ 
	\citep{engwerda2005lq}. Systems that satisfy such conditions are referred to as Euler-Lagrange systems, where the gradient method and conjugate gradient method can be effectively applied to obtain numerical solutions.
	In addition, there exist other approaches, such as the shooting method \citep{hannemann2012verify} which  discretizes the optimal control problems and transform them into nonlinear programming problems.
	
	Nowadays, with the  advancement of computer hardware technology, intelligent computing methods are gradually being integrated into scientific research \citep{sager2009reformulations,koumir2016optimal}.
	In recent years, \citet{raissi2019physics}
	proposed physics-informed neural networks (PINNs). PINNs aim to incorporate physical knowledge  expressed by partial differential equations (PDEs) or ordinary differential equations (ODEs) as prior knowledge into the loss functions of neural networks.
	By penalizing the residuals of PDEs or ODEs, PINNs achieve the numerical solutions of these equations \citep{robinson2022physics,wu2023enhancing,yang2024moving,liu2024diminishing}.
	They can not only solve the forward problems of obtaining numerical solutions for PDEs without data but also address the inverse problems of identifying PDE parameters  in a data-driven manner.
	PINNs  can overcome the limitations of the Runge-Kutta method (numerical solver for ODEs) and the finite element method (numerical solver for PDEs) in handling the large computational cost of high-dimensional problems, effectively alleviating the curse of dimensionality (CoD). The  automatic differentiation technique is utilized by PINNs to directly handle nonlinear problems without relying on any prior assumptions, linearization, or local time-stepping \citep{raissi2019physics}. 
	
	PINNs are also applicable to the study of optimal control problems.
	\cite{hwang2022solving} proposed a 2-phase framework for solving PDE-constrained optimal control problems. 
	\cite{antonelo2024physics}	introduced a PINNs framework for control problems (PINC), which incorporates an autoregressive mechanism and variable time inputs to enable efficient modeling and control of nonlinear dynamic systems.  
	\cite{furfaro2022physics,liu2023pinn,mukherjee2023bridging,shilova2024learning,zhang2024value} presented PINN-based frameworks designed to guarantee solutions to Hamilton-Jacobi-Bellman (HJB) equations within the context of  optimal control problems. In addition to the aforementioned advancements, \citet{schiassi2020physics,schiassi2021class,d2021pontryagin,furfaro2022physics}  explored the application of PINNs for  optimal control problems  of aerospace systems.

	\citet{mowlavi2023optimal} proposed a method for solving PDE-constrained optimal control problems based on PINNs that aims to identify control variables $u(t)$ which minimize a cost objective. The method lies in embedding the appropriate forms of the objective functions $J$ into the loss functions, which enables the minimization of the objective functions to achieve optimality. 
	Nonetheless, constructing an appropriate form of these sub-loss functions may not be easy. 
	Thus, our study, inspired by Euler-Lagrange theorem, proposes Euler-Lagrange theorem combined physics informed neural networks (EL-PINNs) for numerical solutions of Euler-Lagrange systems in optimal control problems. This neural network framework employs three sets of multilayer perceptrons (MLPs) to predict 
	state variables, control variables, and adjoint variables, each taking time 
	$t$  as input. The loss functions of EL-PINNs consist five  key components: 
	the dynamical system equations for state variables and adjoint variables as residual terms; the initial conditions  of   state variables and the terminal conditions of adjoint variables as supervised learning terms; 
	and the residual of  control equations, which penalize to enforce compliance.

	A persistent challenge in conventional PINNs is the loss-imbalance problem: because the composite loss function employs pre-chosen weights, dominant terms often overshadow others, leading to sub-optimal convergence and degraded generalisation.
	PINNs and EL-PINNs can be regarded as a form of multi-task learning, owing to  involving simultaneously optimizing multiple components of the loss function, which is composed of a weighted sum of various constraints.  The performance of multi-task learning heavily depends on the selection of loss weights. However, finding suitable loss weights is challenging, and manually tuning the optimal loss weights is troublesome. \cite{kendall2018multi} addressed multi-task learning in the context of visual scene understanding in computer vision, which involves comprehending both the geometric and semantic information of a scene, and propose an adaptive loss weighting method. Based on the scale of noise in each task, this method introduces the weighting approach grounded in the maximum likelihood estimation of Gaussian noise, which enables an effective balance of weights among tasks, and ultimately achieving superior performance. \cite{xiang2022self,hou2023enhancing} apply this loss weighting method to the task of solving PDEs by PINNs.

	Building on these insights, we propose adaptive weighted Euler-Lagrange theorem combined physics informed neural networks (AW-EL-PINNs), which automatically balances the state, control, and adjoint losses, thereby addressing the loss-imbalance issue inherent in traditional PINN formulations.

	The main contributions of this paper are as follows:
	\begin{enumerate}
		\renewcommand{\labelenumi}{\roman{enumi}.}
		\item We propose EL-PINNs by integrating Euler-Lagrange theorem with  PINNs, in order to achieve  numerical solutions to  Euler-Lagrange systems in optimal control problems unlike previous work in \cite{mowlavi2023optimal} that embeds the corresponding form of the objective functions 
		$J$ directly into the loss functions, we construct multiple ODEs based on Euler-Lagrange theorem as constraints to minimize the loss functions. The essence of our method lies in transforming an optimization problem into a numerical solution problem for ODEs.
		\item We treat the entity of  ODE constraints as a multi-task learning problem. Considering the difficulty of weighting multiple tasks, we incorporate an adaptive loss weighting method and propose AW-EL-PINNs. This
		eliminates the inconvenience and
		achieves superior performance, as demonstrated by the experimental results.
		\item Six numerical examples were solved, including a LQ optimal control problem, a LQ nonzero-sum differential game, and four nonlinear optimal control problems. The results demonstrate that both AW-EL-PINNs and EL-PINNs can achieve highly accurate predicted trajectories. Notably, for nonlinear systems without analytical solutions, the proposed algorithm provides accurate solutions, which offers  potential options for problems with physical applications.
	\end{enumerate}

	The remainder of this paper is structured as follows:
	Section \ref{Sec:03} introduces the fundamental concepts including optimal control problems,  Euler-Lagrange theorem, and multi-task learning paradigm. Section \ref{aw-pinns-el} introduces the detailed architecture of AW-EL-PINNs;
	Section \ref{Sec:04} presents six numerical examples designed to evaluate the performance of AW-EL-PINNs, which demonstrates the  effectiveness through comparisons with other methods;
	Section \ref{Sec:conclusion} presents an integrated analysis of experimental results demonstrating the superior performance of AW-EL-PINNs, concludes the study, and outlines prospective research directions.

	\section{Fundamental Conceptions and Necessary Conditions}\label{Sec:03}
	
	\subsection{Optimal control}
	
	Consider a  system given as: 
	
	\begin{align}\label{sys1}
		\begin{cases}
			\dot{x}(t)=f(t,x(t),u(t)), \quad t\in(t_0,t_f];\\
			x(t_0)=x_0
		\end{cases}
	\end{align}
	with performance indicator
	
	\begin{align}\label{cost1}
		J(t,x(t),u(t))=\int^{t_f}_{t_0}L(t,x(t),u(t))\mathrm{d}t+\Phi(x(t_f)),
	\end{align}  
	where $x:[t_0,t_f]\rightarrow\mathbb{R}^n$ is the state variable of the system;
	$u\in U\subset\mathbb{R}^m$ is the control variable where $U$ is the admissible control set; $f:[t_0,t_f]\times\mathbb{R}^n\times U\rightarrow\mathbb{R}^n $ is the changing rate of state and 
	$x_0$ denotes the initial state of $x(t)$; 
	$J(t,x(t),u(t))$ is the performance indicator with the instantaneous utility function $L(t,x(t),u(t))$ and the terminal utility function $\Phi(x(t_f))$.

	The aim of optimal control problems is to find an optimal control variable $u^*(t)$ which minimizes the performance indicator, i.e., 
	\begin{align}
		J(t,x^*(t),u^*(t))\leq  J(t,x(t),u(t)),
	\end{align}
	where $x^*(t)$ is the optimal state under $u^*(t)$. 
	In other word, the optimal control problems can be described as: 
	\begin{align}
		u^*(t)=\arg\min_{u(\cdot)\in U} J(t,x(t),u(t)).
	\end{align}
	And the minimized performance indicator is called as objective function.

	\begin{assumption}
		The functions $f$, $L$ and $\Phi$ are Lipschitz continuous, i.e., if there exists  constants $\mathcal{G},\mathcal{F}>0$, such that
		\begin{align*}
			&\left\|g(t,x,u)-g(t,y,u)\right\|\leq\mathcal{G}\|x-y\|,\\
			&\left\|\Phi(x)-\Phi(y)\right\|\leq\mathcal{F}\|x-y\|,
		\end{align*}
		where $g$ represents both $f$ and $L$, and $x,y\in\mathbb{R}^n$.
	\end{assumption}
	
	\begin{assumption}
		Functions $f,L,\Phi$ are bounded.
	\end{assumption}
	
	\begin{assumption}
		$f,L$ are  differentiable on $t$ and $x$, and $\Phi$ is  differentiable on $x$.
	\end{assumption}

	\subsection{Pontryagin minimum principle}
	To establish the necessary conditions for optimal control, the Hamiltonian function is  defined formally:
	\begin{align}
		H(t,x(t),u(t),p(t))&:=L(t,x(t),u(t))\nonumber\\
		&\quad+p^T(t)f(t,x(t),u(t)),
	\end{align}
	where $p(t)$ is the adjoint variable as known as Lagrangian multiplier.

	\begin{theorem}\label{theorem of PMP on optimal control} (Pontryagin minimum principle, \citet{engwerda2005lq})
		Let $u^*(t)$ be an optimal control and $x^*(t)$ be the corresponding state. Then, there exists an adjoint variable $p^*(t)\neq 0$ and following equations holding:
		\begin{subequations}
			\begin{align}
				&\dot{x}(t)=f(t,x(t),u(t)), \quad t\in(t_0,t_f];\\
				&x(t_0)=x_0;\\
				&\dot{p}(t)=-\frac{\partial H(t,x(t),u(t),p(t))}{\partial x},\quad t\in[t_0,t_f);\\
				&	p(t_f)=\frac{\partial \Phi(x(t))}{\partial x}\Big|_{t=t_f};\\
				&H(t,x^*(t),u^*(t),p^*(t))\leq H(t,x(t),u(t),p(t)).
			\end{align}
		\end{subequations}
		
	\end{theorem}

	\begin{theorem}\label{Euler-Lagrange}(Euler-Lagrange
		theorem, \citet{engwerda2005lq})
		If the control $u^*(t)$ is the optimal control and $x^*(t)$, $p^*(t)$ are corresponding state and costate, then it is necessary that:
		\begin{subequations}
			\begin{align}
				&\dot{x}(t)=f(t,x(t),u(t)), \quad t\in(t_0,t_f];\\
				&x(t_0)=x_0;\\
				&\dot{p}(t)=-\frac{\partial H(t,x(t),u(t),p(t))}{\partial x},\quad t\in[t_0,t_f);\\
				&	p(t_f)=\frac{\partial \Phi(x(t))}{\partial x}\Big|_{t=t_f};\\
				&\frac{\partial H(t,x^*(t),u^*(t),p^*(t))}{\partial u^*}=0,
			\end{align}
		\end{subequations}
		which reduces the dynamic optimization problem to a static optimization problem.
		
	\end{theorem}

\subsection{Fundamental model of PINNs}

PINNs  are a  method that combines optimization principles with deep learning to obtain numerical solutions for  PDEs or ODEs. Typically, spatio-temporal information serves as inputs to PINNs. Yet, within the framework of this study, we only consider temporal information
$t$ as the sole input to  PINNs.
We consider the fundamental  equations as follows:
\begin{align}\label{pde}
\begin{cases}
r(t,y):=y_t+\mathcal{F}[y]=0, \quad t\in(t_0,t_f)\\
y(t_0)=y_0 \quad \textnormal{or} \quad y(t_f)=y_f,
\end{cases}
\end{align}
where $r(t,y)$ represents the residual function; $y$ denotes the solution of \eqref{pde}; $\mathcal{F}$ is the differential operator; $y_0$ and $y_f$ denote the state value at initial time and terminal time, respectively.

{\color{black} MLP} is the foundational neural network model for this study.
The forward propagation process of an 
$M$-layer neural network  can be presented  as shown below:
\begin{subequations}\begin{align}
&z^0=t;\\
&z^j=\phi\left(w^jz^{j-1}+b^j\right);\\
&y_{NN}=z^M=w^Mz^{M-1}+b^M,
\end{align}\end{subequations}
where $z^j$ represents the output of the neurons in the 
$j$-th layer; $w^j$ and $b^j$ denote the weight matrix and bias vector of the 
$j$-th layer, respectively; $\phi(\cdot)$ is the activation function of the neural network, which is typically chosen as the hyperbolic tangent function ($\tanh$) in PINNs due to its smoothness and differentiability; and $y  _{NN}$ is the predicted output of neural network.

Based on the aforementioned issues, the loss function for PINNs is constructed by mean square error (MSE) as:
\begin{align}\label{loss_original}
\textnormal{Loss}=\omega_rL_r+\omega_iL_i+\omega_fL_f,
\end{align}
where 
\begin{subequations}
\begin{align}
&L_r=\frac{1}{N_r}\sum_{j=1}^{N_r}\Big|\dot{y}^j_{NN}(t_j;\alpha)+\mathcal{F}^j\left[y^j_{NN}(t_j;\alpha)\right]\Big|^2;\label{xresidual}\\
&L_i=\frac{1}{N_i}\sum_{j=1}^{N_i}\Big|
y^j(t_0;\alpha)-y^j_0\Big|^2;\label{xinitial}\\
&L_f=\frac{1}{N_f}\sum_{j=1}^{N_f}\Big|
y^j(t_f;\alpha)-y^j_f\Big|^2,\label{xfinal}
\end{align}
\end{subequations}
and $\omega_r,\omega_i,\omega_f$ are the weights for the respective sub-loss functions. These loss weights are used to balance the various components within the loss function, where    the importance of different optimization components are reflected.

By minimizing the loss function and employing the backpropagation algorithm to optimize the neural network's parameters $\alpha=(w,b)$ with learning rate $\beta_\alpha$, such that:
\begin{align}
\alpha^{k+1}=\alpha^k-\beta_\alpha^k\nabla_\alpha\textnormal{Loss},
\end{align}
the network progressively adjusts the weights and biases during the training process to reduce the discrepancy between the predicted values and the true values. When the loss function stabilizes, the resulting neural network's  final predicted solution represents the numerical solution.

\begin{remark}
Properly configured loss weights can not only effectively enhance the predictive efficiency of the neural networks but also ensure that the model can fully leverage the features of each component when dealing with multi-task or multi-objective optimization. This weight configuration helps optimize the training process, and thereby improves the model's generalization ability and prediction accuracy. Nevertheless, selecting appropriate loss weights is a highly challenging task that often results in troublesome learning costs. To address the difficulties in choosing loss weights, we will introduce an adaptive loss weight update method that ensures the accuracy of model predictions while reducing computational costs.
\end{remark}

\subsection{Multi-task adaptive  loss weights}\label{Multi-task}
Our study will reference the method in \cite{kendall2018multi} to adaptively update the loss weights, with the fundamental principle derived from the theory of maximum likelihood estimation. Let \( y_{NN} \) and $ y $ represent the predicted and true values of the neural network, respectively. We assume that the output of the neural network can be modeled as a Gaussian distribution with a mean of \( y_{NN} \) and a variance of $\sigma^2 $:

\begin{align}
p(y | y_{NN}) = \mathcal{N}(y; y_{NN}, \sigma^2).
\end{align}

Here, $ \sigma$ can be interpreted as the noise scale during the operation of the neural network which indicates the level of uncertainty in the output. Assuming independence among different tasks, the likelihood function for the multi-task output can be expressed as:

\begin{align}\label{likelihood}
p(y_1, y_2, \cdots, y_M | y_{NN}) = \prod_{j=1}^{M} p(y_j | y_{NN}),
\end{align}
where $ y_1, y_2, \cdots,y_M $ represent the model outputs (such as residual terms, initial conditions, terminal conditions).

Consequently, taking negative logarithm of \eqref{likelihood}, we obtain the negative log-likelihood function:
\begin{align}
&-\log p(y_1,y_2,\cdots,y_M|y_{NN})\nonumber\\
&=-\log \sum^M_{j=1}p(y_i|y_{NN})\nonumber\\
&=-\log \sum^M_{j=1}\mathcal{N}(y_j;y_{NN},\sigma^2_j)\nonumber\\
&=-\log \sum^M_{j=1}\left(\frac{1}{\sqrt{2\pi\sigma^2_j}}\exp\left(-\frac{1}{2\sigma^2_j}\frac{1}{N_j}\sum^{ N_j}_{k=1}|y_k-y_{NN}|^2\right)
\right)\nonumber\\
&\propto\sum^M_{j=1}\left(\frac{1}{2\sigma^2_j}\frac{1}{N_j}\sum^{ N_j}_{k=1}\left|y_k-y_{NN}\right|^2+\log\sigma_j
\right).
\end{align}

Thereby, according to the above derivation, the multi-task
weighted loss function is presented as:
\begin{align}
\textnormal{Loss}=\sum^M_{j=1}\left(\frac{1}{2\sigma^2_j}\frac{1}{N_j}\sum^{ N_j}_{k=1}\left|y_k-y_{NN}	\right|^2+\log\sigma_j	\right).
\end{align}
In the above expression, the first term represents the data
fitting term adjusted by $\sigma$, and the second term represents the
regularization term, which constrains the noise parameter $\sigma$ to
prevent it from becoming excessively large.
The model’s loss
function is obtained by minimizing the negative log-likelihood
function, thereby 
each task is adaptively assigned  a weight coefficient. The value of $\sigma$ will be assigned an initial value
before the neural network training, which is typically $1$.
 This setting gives each task an equal initial weight, ensuring that no single term dominates the optimization at the start.
 Both $\sigma$ and $\alpha$ are
trained simultaneously, but based on existing experience, the optimization for the two will utilize different optimizers with
varying learning rates.

During the experimental process, we set $s=\log\sigma^2$. This
formulation allows the loss function to avoid division by zero
and provides greater numerical stability. Additionally, $e^{-s}$ is
mapped to the positive domain, yielding valid values for the loss
weights. Eventually, the loss function in PINNs with adaptive
loss weights can be presented as:
\begin{align}
\textnormal{Loss}=e^{-s_r}L_r+e^{-s_i}L_i+e^{-s_f}L_f+s_r+s_i+s_f,
\end{align}
where $L_r,L_i,L_f$ are presented by \eqref{xresidual}, \eqref{xinitial}, \eqref{xfinal}.
And the update rule for $s$ also follows gradient descent: 
\begin{align}
s^{k+1}=s^k-\gamma_s^k\nabla_s\textnormal{Loss},
\end{align}
which enables the model to dynamically adjust the weights of each task based on their learning progress and difficulty.

Specifically, if the loss $L_j$ for a given task $j$ remains persistently large, the optimizer, in an effort to mitigate the contribution of the term $e^{-s_j}L_j$
to the total loss, tends to increase $s_j$.
An increase in $s_j$ signifies an elevation in the model's estimated uncertainty $\sigma^2_j$
for this task, which in turn reduces the task's effective weight $e^{-s_j}L_j$. This behavior can be interpreted as the model de-prioritizing poorly performing tasks, thereby preventing their substantial losses from disproportionately influencing gradient updates and, consequently, facilitating the learning of other tasks. Conversely, if task $j$ is well-learned (i.e., $L_j$ is small), the optimizer, when balancing the term $e^{-s_j}L_j$
against the regularization term  $s_j$,  permits  $s_j$ to remain at a small value.

\subsection{Existing approach review}\label{exi}
In order to introduce our method, we shall revisit the approach for solving optimal control problems as outlined by   \cite{mowlavi2023optimal}. In this study, the authors constructed two MLPs to
predict the state variable $x(t)$ and the control variable $u(t)$, respectively. Considering the optimization objective of minimizing
the performance indicator $J$, the authors propose adding an external term, $L_J$, to penalize the loss function. The new loss function is constructed as:
\begin{align}
\textnormal{Loss}=\omega_rL_r+\omega_iL_i+\omega_fL_f+\omega_JL_J, 
\end{align}
where the first three terms are derived from \eqref{loss_original}, $L_J$ is selected based on the form of $J$, and $w_J$ is the loss weight of objective term. When the loss
function converges and stabilizes during iterations, the two neural networks
output the optimal control variable along with its corresponding
state trajectory.

Based on the adaptive loss weighting method for multi-tasks described in \ref{Multi-task}, an AW-PINNs approach for solving optimal control problems is proposed with its loss function constructed as:
\begin{align}\label{38}
\textnormal{Loss}&=e^{-s_r}L_r+e^{-s_i}L_i+e^{-s_f}L_f+e^{-s_J}L_J\nonumber\\
&\quad +s_r+s_i+s_f+s_J.
\end{align}

\section{Guideline of   AW-EL-PINNs}\label{aw-pinns-el}

\begin{figure*}[htpb]
\centering
\includegraphics[width=1\linewidth]{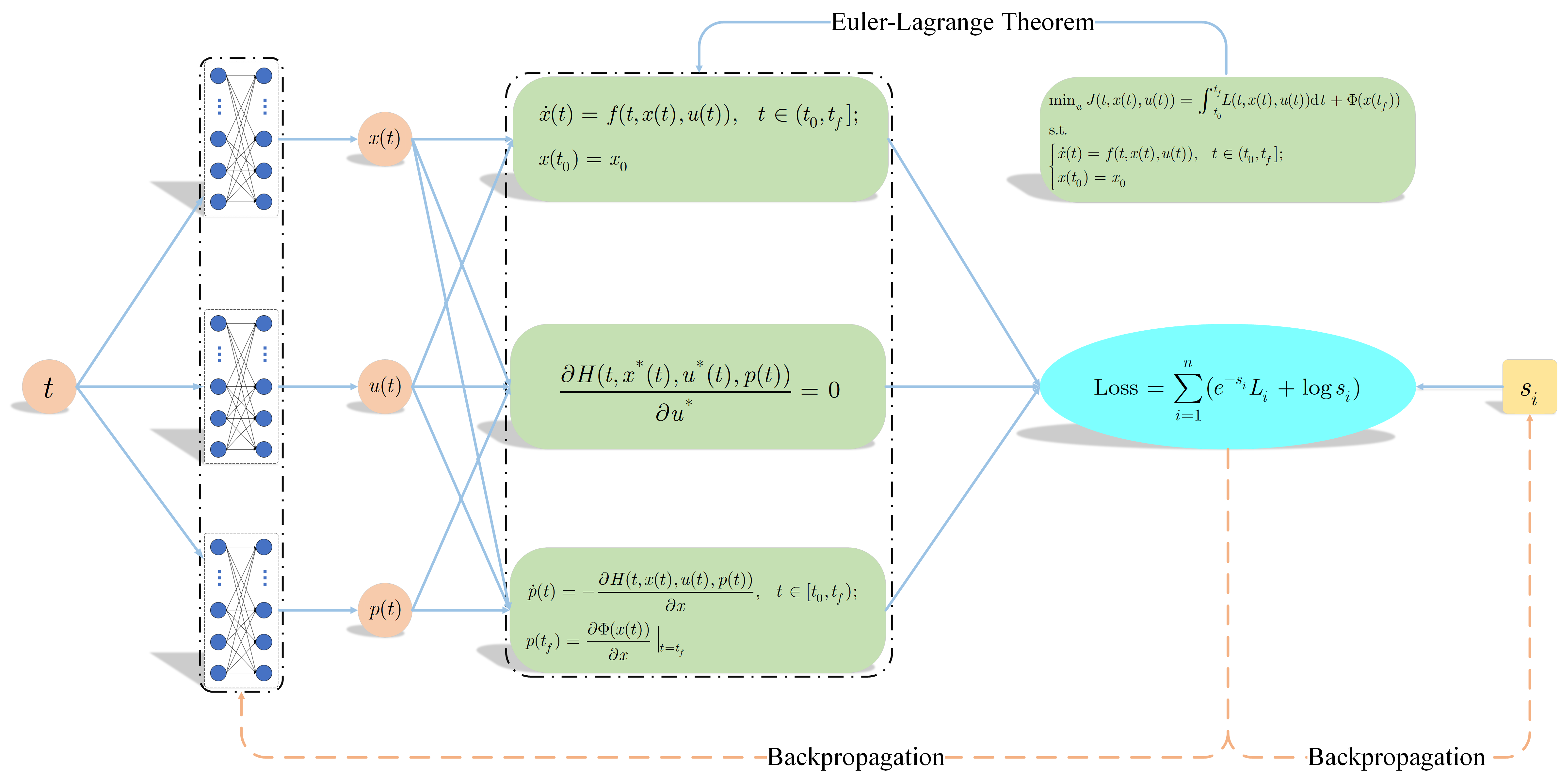}
\caption{Structure of AW-EL-PINNs.}
\label{Structure}
\end{figure*}

Inspired by the above content and Theorem \ref{Euler-Lagrange}, we propose EL-PINNs, a neural network model that incorporates the Euler-Lagrange theorem. This neural network framework consists of three sets of MLPs to predict
$x(t)$, $u(t)$, and $p(t)$, where time $t$ is taken as input for each set. By leveraging the state equations, the adjoint equations, and the control equations, EL-PINNs fit the implicit relationships between these variables using neural networks. The three sets of neural networks are marked as
$NN_x$, $NN_u$, $NN_p$, separately. The loss function of EL-PINNs comprises five components, including the residual of the state functions, along with initial conditions treated as supervised learning terms; the residual of the adjoint functions, with terminal conditions used as supervised learning terms; and the residual of  control functions as a penalty term. The loss function of EL-PINNs is constructed as 
\begin{align} 
\textnormal{Loss}=\omega_xL_x+\omega_{x_0}L_{x_0}+\omega_pL_p+\omega_{p_f}L_{p_f}+\omega_uL_u,
\end{align}
where, 
\begin{subequations}
\begin{align}
L_x &= \frac{1}{N_i} \sum_{j=1}^{N_i} \left| x_{NN}^j(t_j) - f^j(t_j, x_{NN}(t_j), u_{NN}(t_j)) \right|^2,  \\
L_{x_0} &= \frac{1}{N_b} \sum_{j=1}^{N_b} \left| x_{NN}^j(t_0) - x_0 \right|^2,  \\
L_p &= \frac{1}{N_i} \sum_{j=1}^{N_i} \left| p_{NN}^j(t_j) + \frac{\partial H^j}{\partial x_{NN}(t_j)} \right|^2,  \\
L_{p_f} &= \frac{1}{N_b} \sum_{j=1}^{N_b} \left| p_{NN}^j(t_f) - \frac{\partial \Phi(x_{NN}^j(t_f))}{\partial x_{NN}^j(t_f)} \right|^2, \\
L_u &= \frac{1}{N_i} \sum_{j=1}^{N_i} \left| \frac{\partial H^j}{\partial u_{NN}(t_j)} \right|^2,
\end{align}
\end{subequations}
and $\omega_x,\omega_{x_0},\omega_p,\omega_{p_f},\omega_u$ are loss weights.
Utilizing the adaptive loss weighting technique in \ref{Multi-task}, we propose AW-EL-PINNs with the  loss function  constructed as follow, eventually:
\begin{align}\label{sunshi}
\textnormal{Loss}&=e^{-s_x} L_x + e^{-s_{x_0}} L_{x_0} + e^{-s_p} L_p + e^{-s_{p_f}} L_{p_f} + e^{-s_u} L_u \nonumber\\
&\quad+ s_x + s_{x_0} + s_p + s_{p_f} + s_u.
\end{align}
The structure of AW-EL-PINNs is depicted in Fig. \ref{Structure}.

\section{Numerical Illustration and Performance Comparison}\label{Sec:04}
In this section, we present six numerical simulations to validate the effectiveness of   AW-EL-PINNs in solving Euler-Lagrange systems in  optimal control problems, four of which have analytical solutions, while the remaining two lack  analytical expressions.   Within the deep learning framework, we train AW-EL-PINNs and  EL-PINNs in \ref{aw-pinns-el}, AW-PINNs and PINNs in \ref{exi} with the same epochs, deriving the absolute error, relative $\mathcal{L}_2$ error and the variation trajectory of relative $\mathcal{L}_2$ error in order to demonstrate the higher accuracy, better stability and faster convergence of AW-EL-PINNs. The absolute error and relative $\mathcal{L}_2$ error are defined as:

\begin{align}
&\textnormal{Absolute Error}=\left|y_{NN}(t_j)-y(t_j)\right|^2\\
&\textnormal{Relative $\mathcal{L}_2$ Error}=\frac{\sqrt{\sum^M_{j=1}\left|y_{NN}(t_j)-y(t_j)\right|^2}}{\sqrt{\sum^M_{j=1}\left|y(t_j)\right|^2}}
\end{align}
The absolute error reflects the magnitude of the error at a single point but it is less sensitive to the overall distribution of the data. A smaller absolute error indicates that the prediction is closer to the exact values.  The relative $\mathcal{L}_2$ error, on the other hand, reflects the overall error magnitude and takes into account the scale differences in the data, which enables it to provide a fairer assessment when the data spans a wide range. It measures the distance between the predicted solutions and the exact solutions. Therefore, we use the absolute error to evaluate the accuracy of the models and the 
relative $\mathcal{L}_2$ error to assess their stability. Under the adaptive weighting framework, the loss functions include regularization terms, which  result in negative values during the evolution, whereas the loss functions of the baseline model, based on MSE, remain non-negative. Therefore, we do not use changes of loss functions to evaluate model convergence, and on contrary, we assess  based on evolution of  relative $\mathcal{L}_2$ error.

In all subsequent experiments,  we use Pytorch \citep{paszke2019pytorch} as the framework of deep learning.
The same parameter settings
are applied unless otherwise stated.
We employ three sets of MLPs
to predict the state variables $x(t)$, control variables $u(t)$, and
adjoint variables $p(t)$, respectively. Specifically, the network
$NN_x$ is configured with four layers, each containing 50 neurons.
The network $NN_u$ consists of four layers with 30 neurons per
layer. Besides, the network $NN_p$ is designed with four layers, each
comprising 40 neurons. The initial values of $s_k$ are assigned
as random numbers drawn from normal distributions with the
mean of $0$ and the variance of $1$. Uniform sampling of 1000
points is performed in the interval $(t_0,t_f)$, namely $N_i = 1000$.
Additionally, only one point is sampled at both the initial time $t_0$
and the final time $t_f$, namely $N_b = 1$. Optimizer 1 is  chosen to
update the parameters $\alpha$ of the neural networks, with the Adam optimizer
and the learning rate of $1e-4$. Optimizer 2, which is used to
update $s_i$, also employs the Adam optimizer, with a learning
rate of  $1e-3$. 
In the comparative experiments, the configurations of the three traditional algorithms are as follows: For the gradient method and the conjugate gradient method, the time step is set to $1e-3$, the maximum number of iterations is $500$, the convergence tolerance is 
$1e-3$, and the initial value is set as a zero vector.
For the shooting method, the number of shooting points is $20$, the maximum number of function evaluations is 
$1e+5$
, the convergence tolerance is 
$1e-6$, and the initial guess is also set as a zero vector.

\begin{example}\label{Linear quadratic (LQ) optimal control problem} (LQ optimal control problem \citep{hager1990multiplier})

Consider  the following problem:
\begin{subequations}
\begin{align}
	&\min_uJ(x(t),u(t))=\frac{1}{4}\int^1_0 1.25x(t)^2+x(t)u(t)+u(t)^2\mathrm{d}t,\\
	&\textnormal{s.t.}\nonumber\\
	&\dot{x}(t)=0.5x(t)+u(t),\\
	&x(0)=1,
\end{align}
\end{subequations}	
with the analytical solutions as:
\begin{subequations}
\begin{align}
	&x^*(t)=\frac{\cosh(1-t)}{\cosh(1)},\\
	&u^*(t)=-x^*(t)\left[\tanh(1-t)+0.5\right].
\end{align}
\end{subequations}

Initially, the presented AW-EL-PINNs are utilized for solving.   Following the setup of loss function in section \ref{aw-pinns-el}, after 30000 iterations, the trajectories of state variable $x(t)$ and control variable $u(t)$ are shown in Fig. \ref{fig_ex1:aw_pmpinns}. In these figures, the 
solid lines represent the predicted solutions by the networks, while the dashed ones represent the  analytical solution trajectories.
And the updates of loss weights in AW-EL-PINNs are shown in Fig. \ref{fig_ex1:exp_aw_pmpinn}. 
\begin{figure}[H]
	\centering
	\begin{minipage}{0.48\textwidth}
		\centering
		\includegraphics[width=0.48\textwidth]{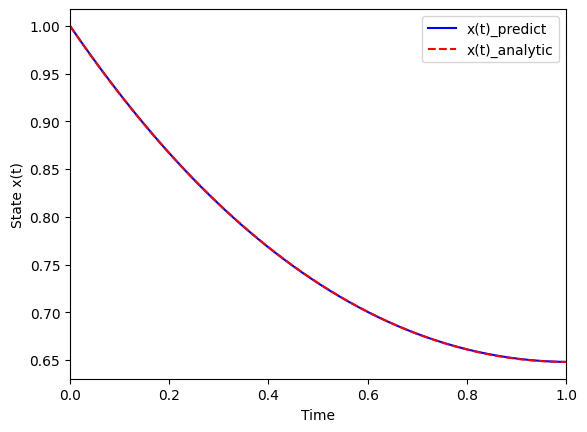}
		\hfil
		\includegraphics[width=0.48\textwidth]{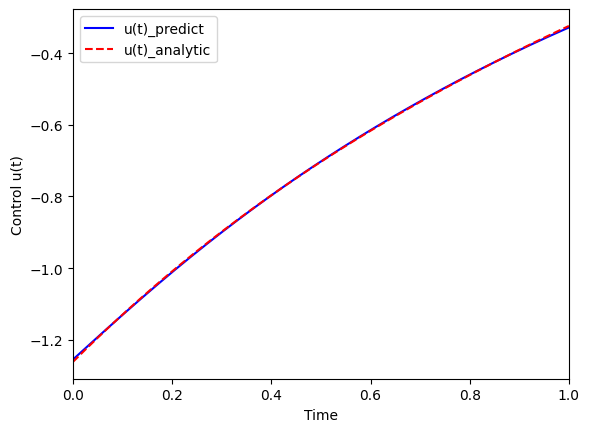}
		\caption{$x(t)$ and $u(t)$ in Example \ref{Linear quadratic (LQ) optimal control problem} solved by AW-EL-PINNs.}
		\label{fig_ex1:aw_pmpinns}
	\end{minipage}
\end{figure}

\begin{figure}[H]
	\centering
	\includegraphics[width=0.85\linewidth]{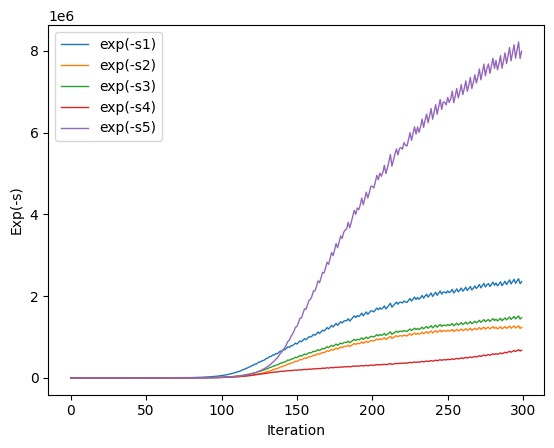}
	\caption{$e^{-s_k}$ of AW-EL-PINNs in Example \ref{Linear quadratic (LQ) optimal control problem}.}
	\label{fig_ex1:exp_aw_pmpinn}
\end{figure}

Subsequently, we consider a set of fixed loss weights $\omega_1=\cdots=\omega_5=1$, with, consequently, the trajectory plots given by Fig. \ref{fig_ex1:pmpinns}.
\begin{figure}[H]
	\begin{minipage}{0.48\textwidth}
		\centering
		\includegraphics[width=0.48\textwidth]{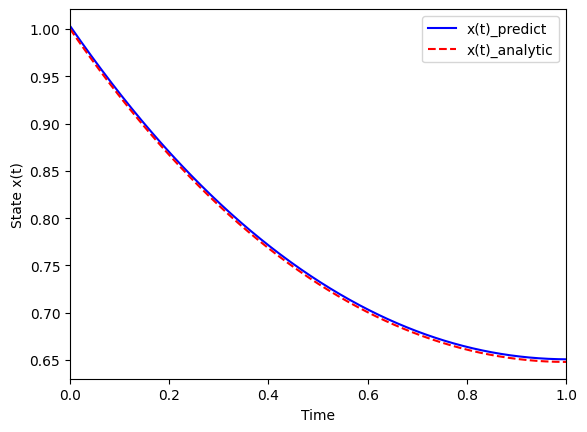}
		\hfil
		\includegraphics[width=0.48\textwidth]{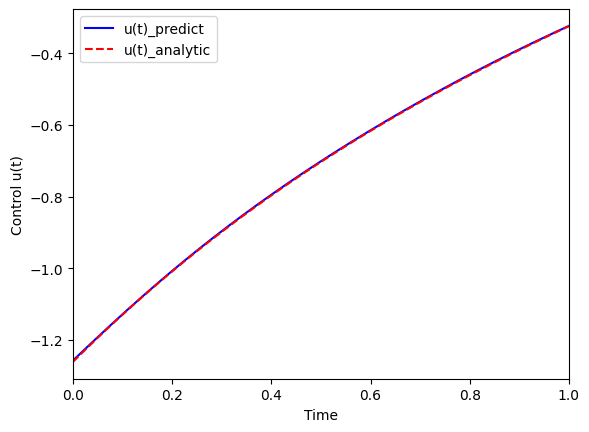}
		\caption{$x(t)$ and $u(t)$ in Example \ref{Linear quadratic (LQ) optimal control problem} solved by EL-PINNs.}
		\label{fig_ex1:pmpinns}
	\end{minipage}
\end{figure}

In succession, based on the model in \cite{mowlavi2023optimal} and incorporated with adaptive loss weights, we design AW-PINNs to solve this problem. We set $L_J$ and the loss function  as:
\begin{align}
	&L_J=\left(\frac{1}{N_i}\sum^{N_i}_{i=1}\frac{1}{4} 1.25x^j_{NN}(t)^2+x^j_{NN}(t)u^j_{NN}(t)+u^j_{NN}(t)^2\right)^2,
\end{align}	
and the loss function is set up as \eqref{38}. 
Similarly, the fixed loss weights are chosen as $\omega_1=\omega_2=1$, $\omega_J=0.05$. We get the trajectory plots in Fig. \ref{fig_ex1:aw_pinns}-Fig. \ref{fig_ex1:pinns}.

\begin{figure}[H]
	\begin{minipage}{0.48\textwidth}
		\centering
		\includegraphics[width=0.48\textwidth]{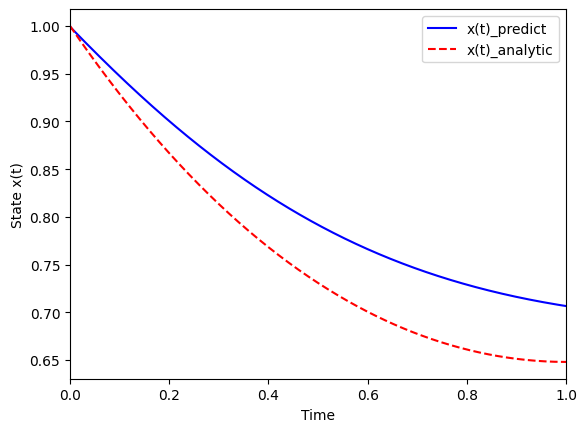}
		\hfil
		\includegraphics[width=0.48\textwidth]{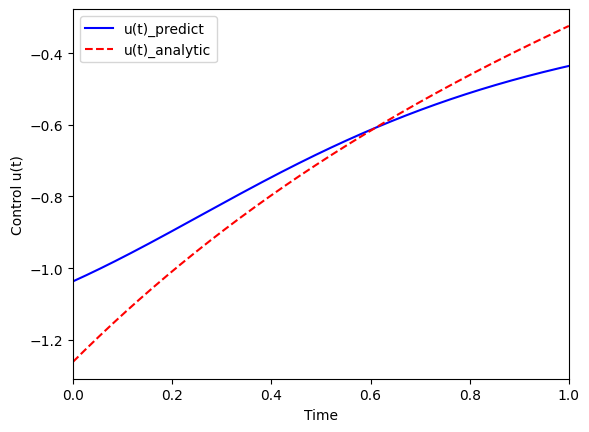}
		\caption{$x(t)$ and $u(t)$ in Example \ref{Linear quadratic (LQ) optimal control problem} solved by AW-PINNs.}
		\label{fig_ex1:aw_pinns}
	\end{minipage}
\end{figure}
\begin{figure}[H]
	\begin{minipage}{0.48\textwidth}
		\centering
		\includegraphics[width=0.48\textwidth]{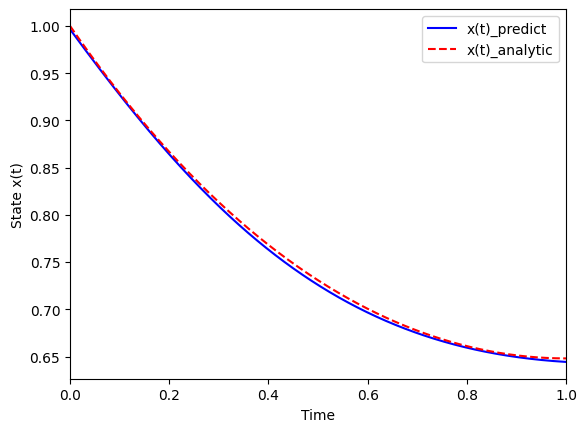}
		\hfil
		\includegraphics[width=0.48\textwidth]{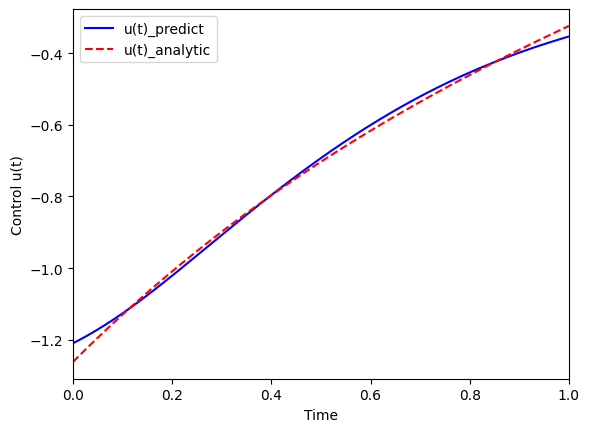}
		\caption{$x(t)$ and $u(t)$ in Example \ref{Linear quadratic (LQ) optimal control problem} solved by PINNs.}
		\label{fig_ex1:pinns}
	\end{minipage}
\end{figure}

Except for  Fig. \ref{fig_ex1:aw_pinns}, the predicted trajectories by neural networks in the other three groups closely coincide with the analytical solution trajectories, with AW-EL-PINNs achieving higher accuracy in its predictions. 

To further compare the accuracy of the four models, we present the box plots of absolute errors  depicted  in Fig. \ref{fig_ex1:box}:
\begin{figure}[H]
	\centering
	\includegraphics[width=0.85\linewidth]{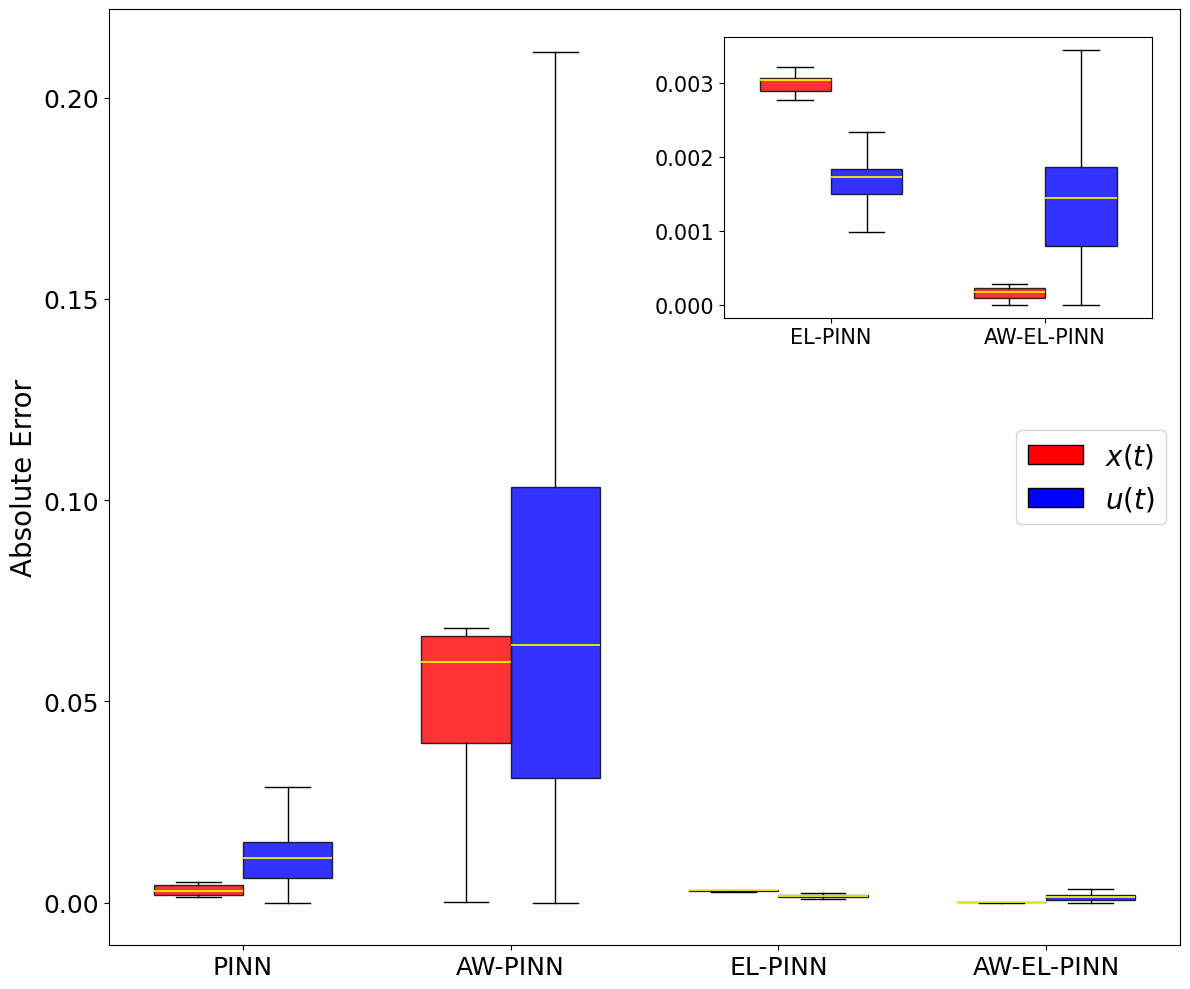}
	\caption{Box plots of absolute errors in Example \ref{Linear quadratic (LQ) optimal control problem}.}
	\label{fig_ex1:box}
\end{figure}

From Fig. \ref{fig_ex1:box}, it can be analyzed that AW-EL-PINNs and EL-PINNs have smaller maximum absolute errors compared to AW-PINNs and PINNs.  Furthermore, between AW-EL-PINNs and EL-PINNs, AW-EL-PINNs exhibit smaller maximum and minimum absolute errors for the prediction of $x(t)$, while for the prediction of $u(t)$,  AW-EL-PINNs show larger maximum absolute errors but smaller minimum absolute errors.  It indicates that solving optimal control problems using AW-EL-PINNs and EL-PINNs provides higher accuracy compared to AW-PINNs and PINNs.  For EL-PINNs with appropriately selected loss weights, they may achieve higher accuracy than AW-EL-PINNs in predicting certain variables.

In addition,  we provide the bar charts of relative $\mathcal{L}_2$ errors in Fig. \ref{fig_ex1:bar} to demonstrate the stability of the solutions obtained by the four models.

\begin{figure}[H]
	\centering
	\includegraphics[width=0.85\linewidth]{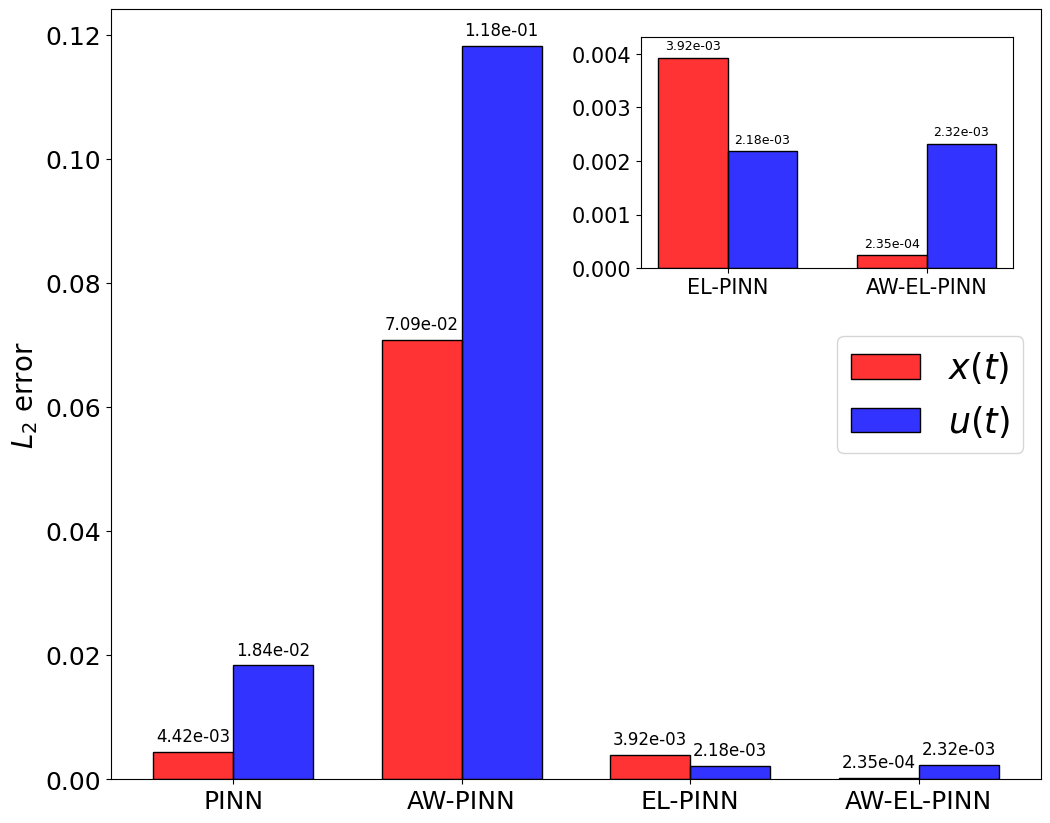}
	\caption{Bar charts of 
		relative $\mathcal{L}_2$ errors in Example \ref{Linear quadratic (LQ) optimal control problem}.}
	\label{fig_ex1:bar}
\end{figure}
Depicted in Fig. \ref{fig_ex1:bar}, it can be observed that the relative $\mathcal{L}_2$ errors of AW-EL-PINNs and EL-PINNs are significantly smaller than those of AW-PINNs and PINNs. AW-EL-PINNs exhibit smaller 
relative $\mathcal{L}_2$ errors for the prediction of 
$x(t)$ compared to EL-PINNs, but slightly larger relative $\mathcal{L}_2$ errors for the prediction of 
$u(t)$. This indicates that AW-EL-PINNs and EL-PINNs demonstrate higher stability than AW-PINNs and PINNs. Moreover, AW-EL-PINNs are more stable than EL-PINNs with fixed loss weights in prediction.

Eventually, we present the evolution plots of relative $\mathcal{L}_2$ errors during the iterative process to demonstrate the convergence speed of the four models.
\begin{figure}[H]
	
	\centering
	\includegraphics[width=0.85\linewidth]{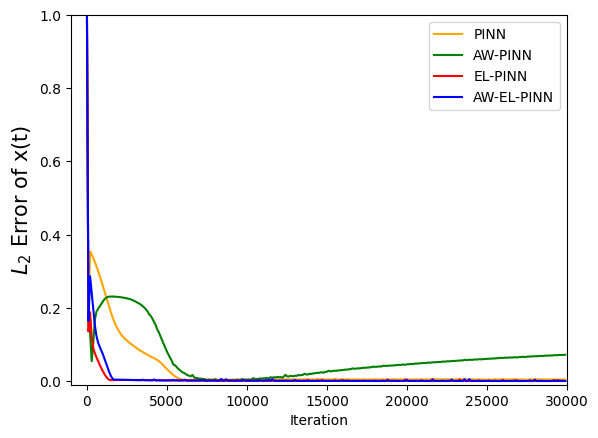}
	\caption{Evolution of $x(t)$'s relative $\mathcal{L}_2$ errors in Example \ref{Linear quadratic (LQ) optimal control problem}.}
	\label{fig_ex1:l2x}
	
	\begin{figure}[H]
	\end{figure}
	
	\centering
	\includegraphics[width=0.85\linewidth]{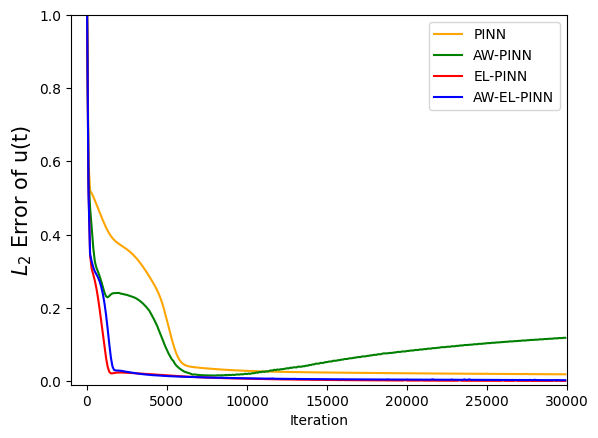}
	\caption{Evolution of $u(t)$'s relative $\mathcal{L}_2$ errors in Example \ref{Linear quadratic (LQ) optimal control problem}.}
	\label{fig_ex1:l2u}
	
\end{figure}

Fig. \ref{fig_ex1:l2x} and Fig. \ref{fig_ex1:l2u} show that both AW-EL-PINNs and EL-PINNs achieve faster convergence, with EL-PINNs exhibiting a faster convergence speed compared to AW-EL-PINNs.

Furthermore, we also consider numerical methods like gradient  method, conjugate gradient method, shooting method, in order to demonstrate the outperform of our method.
Since the maximum absolute error (MAE) and relative $\mathcal{L}_2$ error (RLE) of $x(t)$ and $u(t)$ are the same order of magnitude, which means taking the average values of them does not allow any single variable to dominate the overall metric, the average values between $x(t)$ and $u(t)$ of MAE and RLE for  each method are chosen to be compared and 
shown in  table \ref{tabopc1duibi}.

\begin{table}[h]  
	\centering  
	\begin{tabular}{ccc}  
		\hline
		\hline
		Method& MAE  & RLE   \\
		\hline
		Gradient & 5.96e-02 & 3.11e-02 \\
		Conjugate Gradient & 6.04e-02 & 3.46e-02\\
		Shooting & 2.28e-02 & 1.16e-02 \\
		PINNs & 2.86e-02 & 1.14e-02\\
		AW-PINNs & 1.46e-01 &9.45e-02\\
		EL-PINNs & 3.50e-03 & 3.05e-03\\
		AW-EL-PINNs & {\bf 3.42e-03} & {\bf 1.28e-03}\\
		\hline\hline
	\end{tabular}
	\caption{Maximum absolute error and relative $\mathcal{L}_2$ error of numerical methods in Example \ref{Linear quadratic (LQ) optimal control problem}}  
	\label{tabopc1duibi}  
\end{table}
From the results in the table, it is evident that, among the six numerical methods, our proposed PINN-based approach exhibits superior accuracy, with AW-EL-PINNs achieving the highest precision.

\begin{remark}
	As shown in Fig. \ref{fig_ex1:l2x} and Fig. \ref{fig_ex1:l2u}, AW-PINNs become unstable after approximately 7000 iterations.
	To address this issue, we limited the number of iterations to 7000, and subsequently obtained the corresponding prediction plots for  $x(t)$ and $u(t)$, along with the evolution plots of relative $\mathcal{L}_2$ errors.
	
	\begin{figure}[H]
		\centering
		\begin{minipage}{0.48\textwidth}
			\centering
			\includegraphics[width=0.48\textwidth]{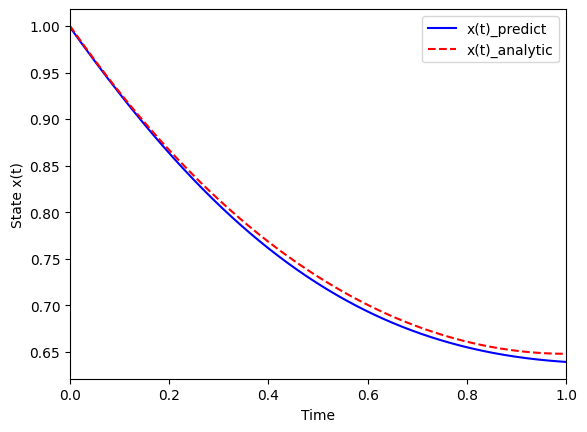}
			\hfil
			\includegraphics[width=0.48\textwidth]{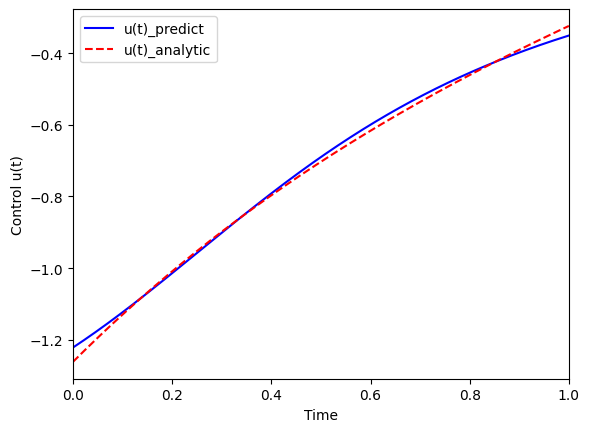}
			\caption{$x(t)$ and $u(t)$ in Example \ref{Linear quadratic (LQ) optimal control problem} solved by AW-PINNS in 7000 iterations.}
			\label{fig_ex1:aw_pinns7000}
		\end{minipage}
	\end{figure}

	It can be observed that after adjusting the number of iterations for AW-PINNs, the predicted solution and the analytical solution trajectories generally align. However, the accuracy is still inferior compared to AW-EL-PINNs.
\end{remark}
\end{example}

	\begin{example}\label{LQ differential game} (LQ open-loop Nash differential game \citep{nikooeinejad2016numerical})
		
		An LQ nonzero-sum differential game with two
		players gives the performance indicators as:
		
		\begin{subequations}
				\begin{align}
						&\min_{u_1}J_1(x(t),u_1(t))=\int^3_0\left(x(t)^2+u_1(t)^2\right)\mathrm{d}t,\\
						&\min_{u_2}J_2(x(t),u_2(t))=\int^3_0\left(4x(t)^2+u_2(t)^2\right)\mathrm{d}t+5x\left(T\right)^2\\
						&\textnormal{s.t.}\nonumber\\
						&\dot{x}(t)=2x(t)+u_1(t)+u_2(t),\\
						&x(0)=1.
					\end{align}
			\end{subequations}
		The analytical solutions of open-loop Nash equilibrium for this problem state as:
		\begin{subequations}
				\begin{align}
						&x^*(t)=e^{-3t},\\
						&u^*_1(t)=-e^{-3t}+e^{-(2t+3)},\\
						&u^*_2(t)=-4e^{-3t}-e^{-(2t+3)}.
					\end{align}
			\end{subequations}

	Initially, we consider solving with AW-EL-PINNs and EL-PINNs. The loss weights for EL-PINNs are set as \( \omega_1 = \cdots = \omega_5 = 1 \), \( \omega_6 = 11.1 \), and \( \omega_7 = 10.1 \). After 60000 iterations, two sets of the trajectory plots are obtained in Fig. \ref{fig_exnash:aw_pmpinns} and Fig. \ref{fig_exnash:pmpinns}.

		\begin{figure}[H]
				\centering
				\begin{minipage}{0.48\textwidth}
						\centering
						\includegraphics[width=0.48\textwidth]{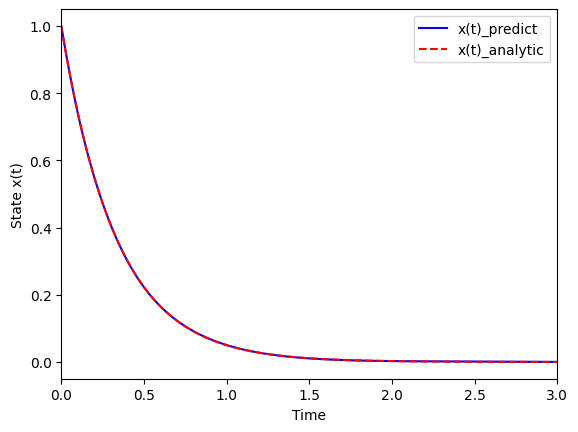}
						\hfil
						\includegraphics[width=0.48\textwidth]{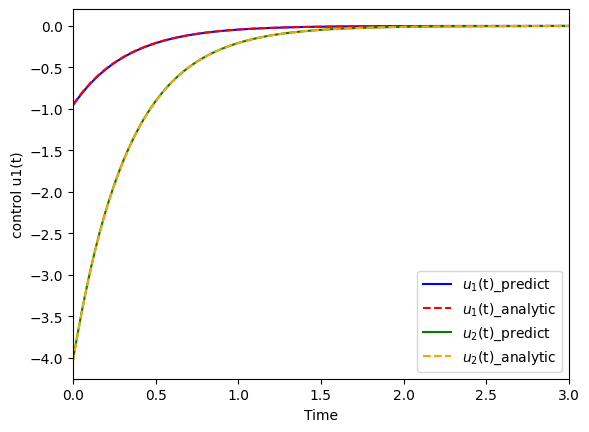}
						\caption{$x(t)$ and $u(t)$ in Example \ref{LQ differential game} solved by AW-EL-PINNs.}
						\label{fig_exnash:aw_pmpinns}
					\end{minipage}
			\end{figure}
		\begin{figure}[H]
				\begin{minipage}{0.48\textwidth}
						\centering
						\includegraphics[width=0.48\textwidth]{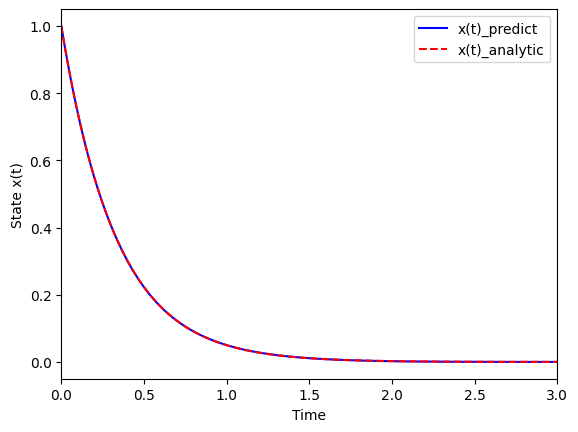}
						\hfil
						\includegraphics[width=0.48\textwidth]{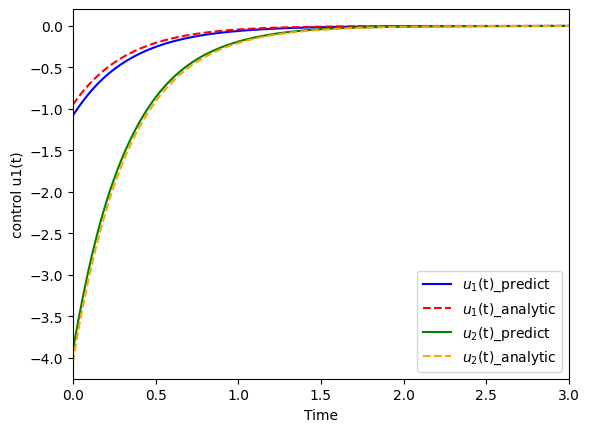}
						\caption{$x(t)$ and $u(t)$ in Example \ref{LQ differential game} solved by EL-PINNs.}
						\label{fig_exnash:pmpinns}
					\end{minipage}
			\end{figure}
			The updates of loss weights of  AW-EL-PINNs and AW-PINNs given in Fig. \ref{fig_exnash:exp_aw_pmpinn}:
		\begin{figure}[H]
				\centering
				\includegraphics[width=0.85\linewidth]{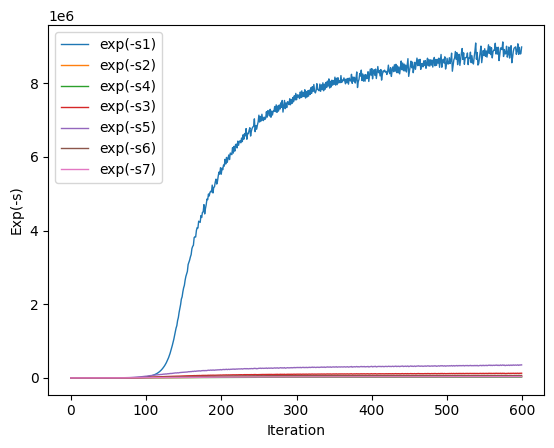}
				\caption{$e^{-s_k}$ of AW-EL-PINNs in Example \ref{LQ differential game}.}
				\label{fig_exnash:exp_aw_pmpinn}
				
			\end{figure}
		Similar to Example \ref{Linear quadratic (LQ) optimal control problem}, we solve this differential game by AW-PINNs and PINNs in \cite{mowlavi2023optimal}. Set $\omega_{J1}=\omega_{J2}=1$, and 
		\begin{subequations}
				\begin{align}
					&L_{J1}=\left(\frac{1}{N_i}\sum^{N_i}_{j=1}\left(x^j_{NN}(t_j)^2+u^j_{1NN}(t_j)^2\right)\right)^2,\\
					&L_{J2}=\left(
					\frac{1}{N_i}\sum^{N_i}_{j=1}\left(4x^j_{NN}(t_j)^2+u^j_{2NN}(t_j)^2\right)+\frac{1}{N_b}\sum^{N_b}_{j=1}5x^j_{NN}(\mathbf{3}[j])
					\right)^2.
				    \end{align}
				\end{subequations}
		Then the trajectories are plotted in Fig. \ref{fig_exnash:aw_pinns} and Fig. \ref{fig_exnash:pinns}.
		\begin{figure}[H]
				\begin{minipage}{0.48\textwidth}
						\centering
						\includegraphics[width=0.48\textwidth]{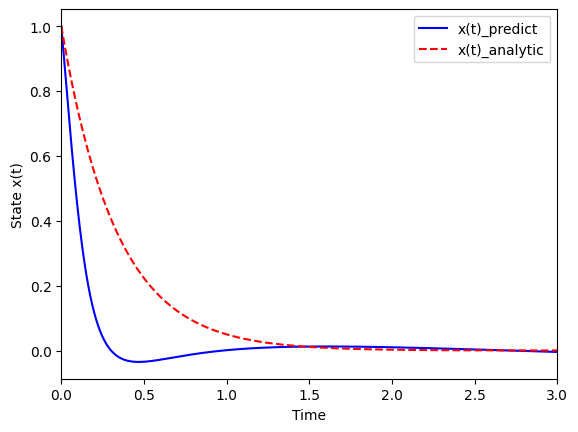}
						\hfil
						\includegraphics[width=0.48\textwidth]{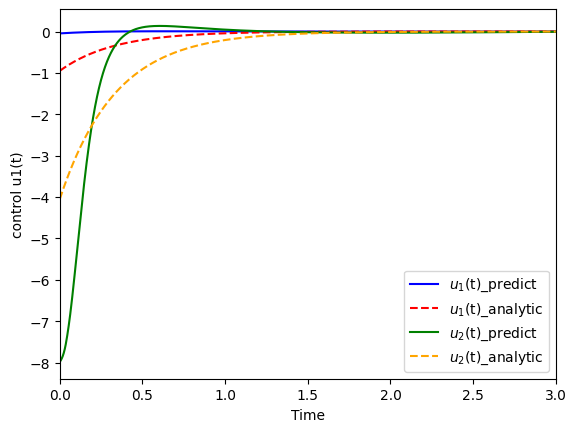}
						\caption{$x(t)$ and $u(t)$ in Example \ref{LQ differential game} solved by AW-PINNs.}
						\label{fig_exnash:aw_pinns}
					\end{minipage}
			\end{figure}
		\begin{figure}[H]
				\begin{minipage}{0.48\textwidth}
						\centering
						\includegraphics[width=0.48\textwidth]{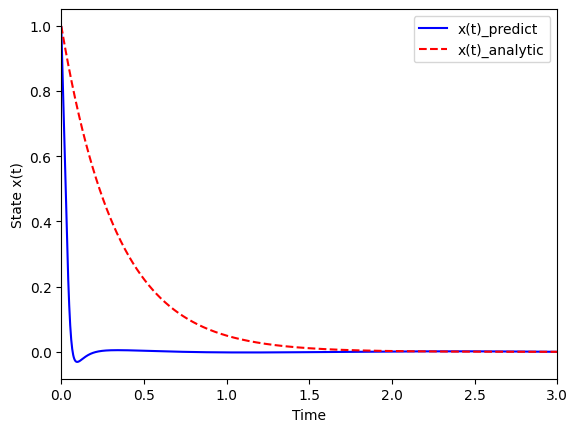}
						\hfil
						\includegraphics[width=0.48\textwidth]{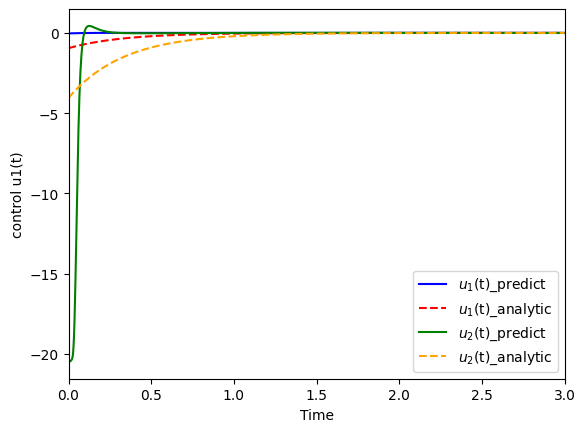}
						\caption{$x(t)$ and $u(t)$ in Example \ref{LQ differential game} solved by PINNs.}
						\label{fig_exnash:pinns}
					\end{minipage}
			\end{figure}

		Based on Fig. \ref{fig_exnash:aw_pmpinns}, Fig. \ref{fig_exnash:pmpinns},  Fig. \ref{fig_exnash:aw_pinns} and Fig. \ref{fig_exnash:pinns}, it can be concluded that AW-PINNs and PINNs fail to solve this problem, while AW-EL-PINNs and EL-PINNs provide highly accurate predicted trajectories.
		To further investigate the accuracy, we also provide the box plots of absolute errors for this problem in Fig. \ref{fig_exnash:box}.
	
	\begin{figure}[H]
				\centering
				\includegraphics[width=0.85\linewidth]{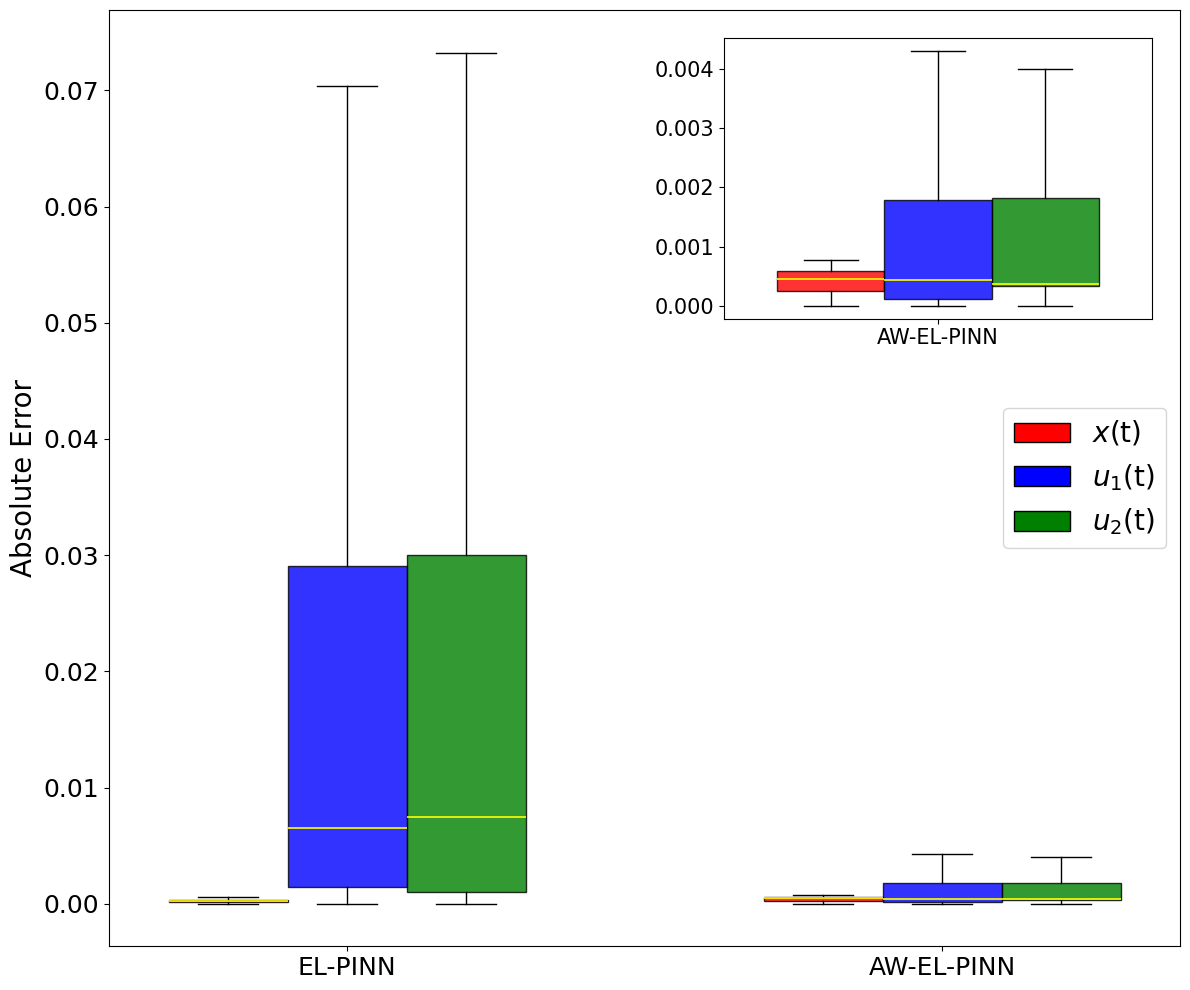}
				\caption{Box plots of absolute errors in Example \ref{LQ differential game}.}
				\label{fig_exnash:box}
				\end{figure}
			
			Fig. \ref{fig_exnash:box} demonstrates  that  AW-EL-PINNs and EL-PINNs achieve comparable accuracy in predicting $x(t)$'s absolute error, while AW-EL-PINNs outperform in predicting $u_1(t)$ and $u_2(t)$.
			In addition to the box plots of absolute errors, we also plot bar charts of relative $\mathcal{L}_2$ errors in Fig. \ref{fig_exnash:bar} to further explore the stability of the models.

		\begin{figure}[H]
			\centering
			\includegraphics[width=0.85\linewidth]{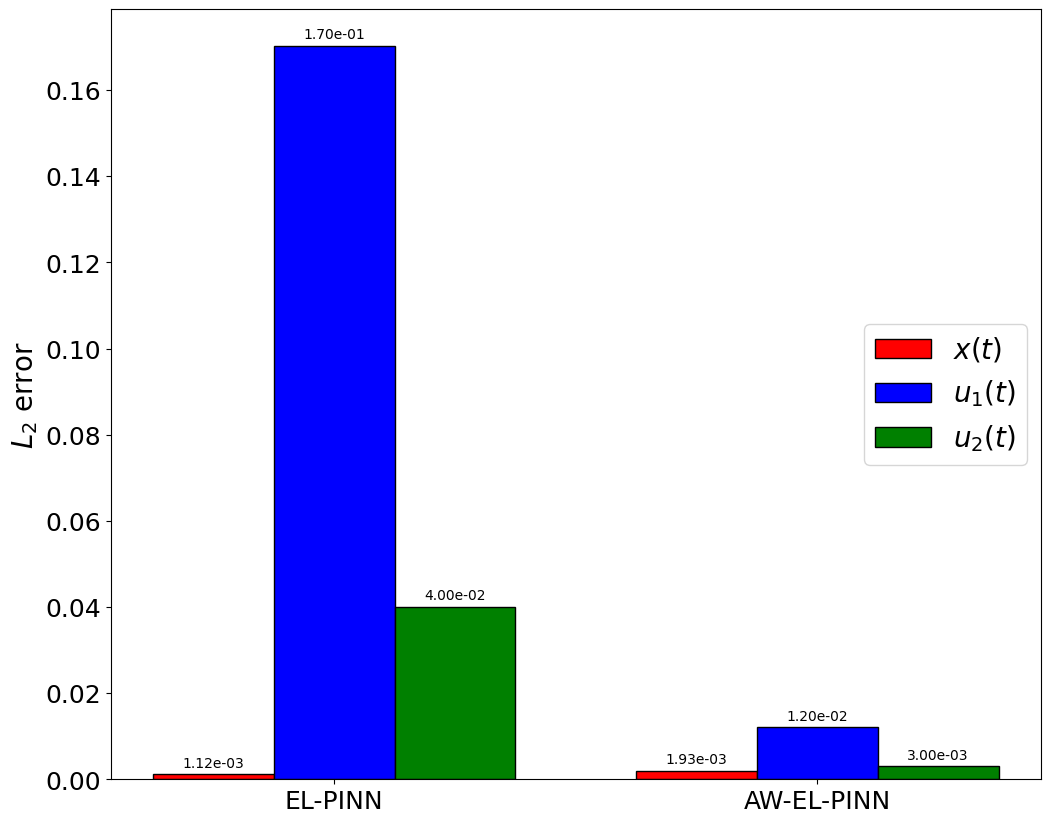}
			\caption{Bar charts of relative $\mathcal{L}_2$ errors in Example \ref{LQ differential game}.}
			\label{fig_exnash:bar}
		\end{figure}
	
	Fig. \ref{fig_exnash:bar} demonstrates that AW-EL-PINNs and EL-PINNs exhibit strong stability. However,  EL-PINNs achieve slightly smaller relative $\mathcal{L}_2$ error in predicting $x(t)$, while AW-EL-PINNs show approximately 10 times the relative $\mathcal{L}_2$ errors of EL-PINNs in predicting $u_1(t)$ and $u_2(t)$.
		
		Finally, we provide the evolution plots of 
		relative $\mathcal{L}_2$ errors throughout the iterative process to illustrate the convergence speeds of AW-EL-PINNs and EL-PINNs.
	
			\begin{figure}[H]
	
					\centering
					\includegraphics[width=0.85\linewidth]{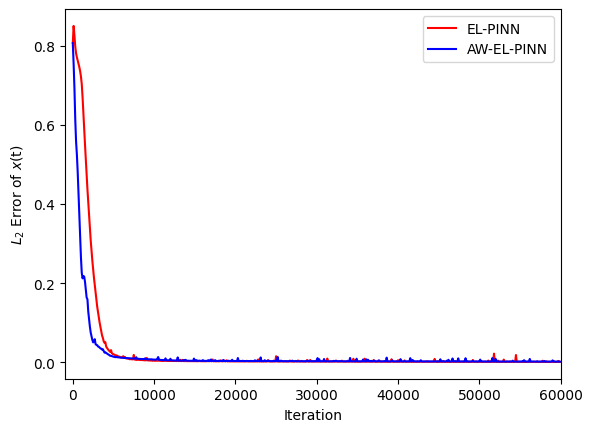}
					\caption{Evolution of $x(t)$'s relative $\mathcal{L}_2$ errors in Example \ref{LQ differential game}.}
					\label{fig_exnash:xL2}
				\end{figure}
		\begin{figure}[H]
					\centering
					\includegraphics[width=0.85\linewidth]{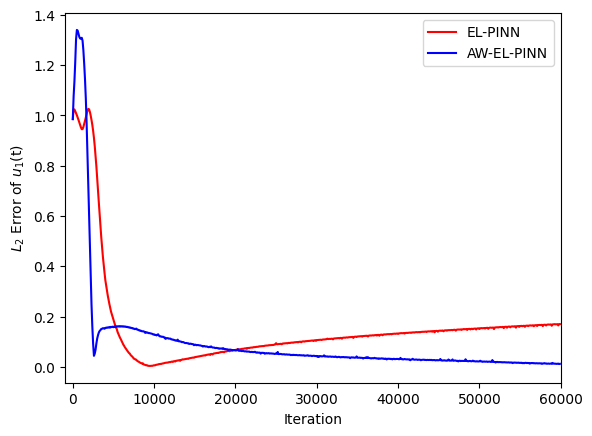}
					\caption{Evolution of $u_1(t)$'s relative $\mathcal{L}_2$ errors in Example \ref{LQ differential game}.}
					\label{fig_exnash:u1L2}
		\end{figure}
\begin{figure}[H]
				\centering
				\includegraphics[width=0.85\linewidth]{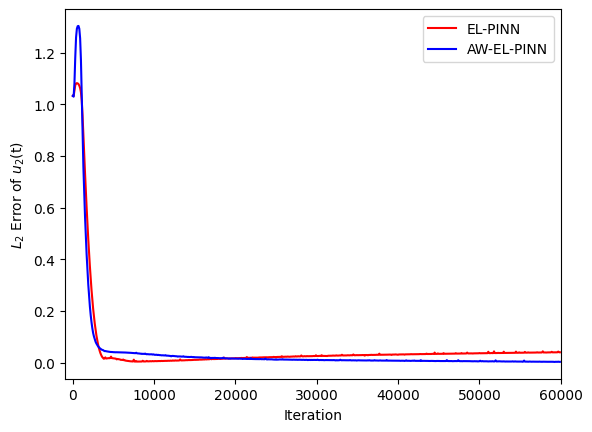}
				\caption{Evolution of $u_2(t)$'s relative $\mathcal{L}_2$ errors in Example \ref{LQ differential game}.}
				\label{fig_exnash:u2L2}
		\end{figure}
		
		Fig. \ref{fig_exnash:xL2} clearly illustrate that AW-EL-PINNs exhibit faster convergence when predicting $x(t)$. In Fig. \ref{fig_exnash:u1L2} and Fig. \ref{fig_exnash:u2L2}, the relative $\mathcal{L}_2$ error curves of EL-PINNs decrease to near 0 more quickly than that of AW-EL-PINNs but then gradually begin to increase.

		\begin{remark}
			According to Fig. \ref{fig_exnash:u1L2}, we observe that after approximately 10000 iterations, the relative $\mathcal{L}_2$ error of $u_1(t)$ starts to gradually increase. Therefore, we attempt to limit the number of iterations for EL-PINNs to 10000, leading to the following trajectory and relative $\mathcal{L}_2$ error plots. 
		\end{remark}
		
		\begin{figure}[H]
			\centering
			\begin{minipage}{0.48\textwidth}
					\centering
					\includegraphics[width=0.48\textwidth]{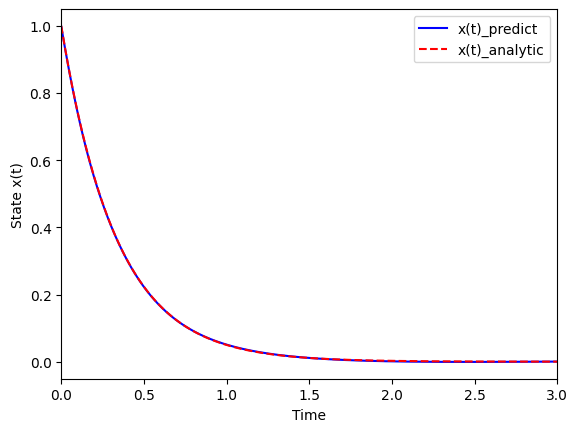}
					\hfil
					\includegraphics[width=0.48\textwidth]{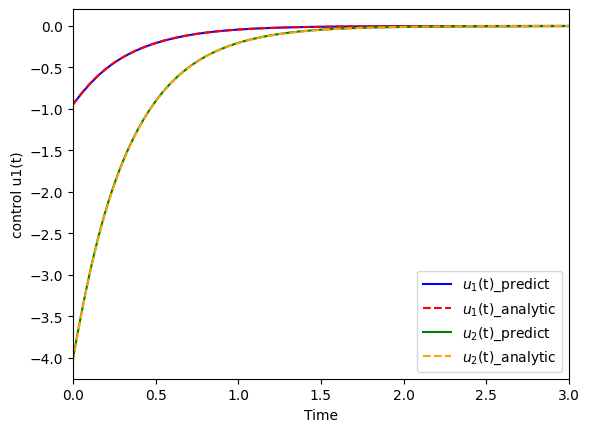}
					\caption{$x(t)$ and $u(t)$ in Example \ref{LQ differential game} solved by EL-PINNs in 10000 iterations.}
					\label{fig_exnash:pmpinns10000}
				\end{minipage}
		\end{figure}
		\begin{figure}[H]
			\centering
			\includegraphics[width=0.85\linewidth]{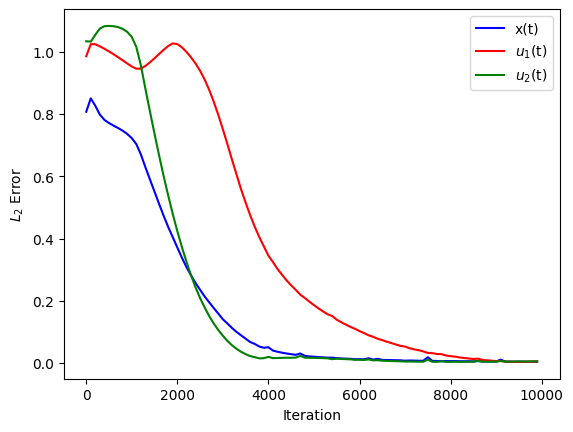}
			\caption{relative $\mathcal{L}_2$ errors in Example \ref{LQ differential game} solved by EL-PINNs in 10000 iterations.}
			\label{fig_exnash:l2_pmpinns7000}
		\end{figure}

		It is evident that for EL-PINNs, selecting an appropriate number of iterations can also yield highly accurate predicted solutions.

	\end{example}

\begin{example}\label{example of nonlinear} (Nonlinear maximum terminal state problem)

Given  the below optimal control problem for $t\in[0,5]$, the objective is to maximize the terminal state:
\begin{subequations}
\begin{align}
	&\max_u J(x(t))=x(t_f),\\
	&\textnormal{s.t.}\nonumber\\
	&\dot{x}(t)=-x(t)+u(t)x(t)-u(t)^2,\\
	&x(0)=1.
\end{align}
\end{subequations}
The corresponding analytical solution is derived as:
\begin{subequations}
\begin{align}
	&x^*(t)=\frac{4}{1+3e^t},\\
	&x^*(t)=\frac{2}{1+3e^t}.
\end{align}
\end{subequations}

The problem aims to maximize the terminal state, which is equivalent to minimizing its negative value. We attempted to solve this problem using both AW-PINNs and PINNs, but neither yielded satisfied results. The primary issue lies in the difficulty of selecting $L_J$. When $-x(t_f)$ is directly added as $L_J$ in the loss function, the model becomes unstable. We also experimented with a series of monotonically increasing and bounded functions defined over the entire real domain as candidates for $L_J$, such as $L_J=\arctan(-x(t_f))$, but the model still exhibited instability. Similarly, when setting $L_J=\textnormal{sigmoid}(-x(t_f))$, the model was stable but failed to produce predictions that aligned with the analytical solution trajectory. Consequently, AW-PINNs and PINNs are unsuitable for solving this problem effectively. In the following discussion, we focus exclusively on the performance of AW-EL-PINNs and EL-PINNs for this problem.
The loss weights of EL-PINNs are setting as $\omega_1=\omega_2=\omega_3=\omega_5=1$, $\omega_4=25$. After 20000 iterations, the predicted solution trajectories of the two models are  visualized as Fig. \ref{fig_exnon:aw_pmpinns}-Fig. \ref{fig_exnon:pmpinns}. 
\begin{figure}[H]
\centering
\begin{minipage}{0.48\textwidth}
	\centering
	\includegraphics[width=0.48\textwidth]{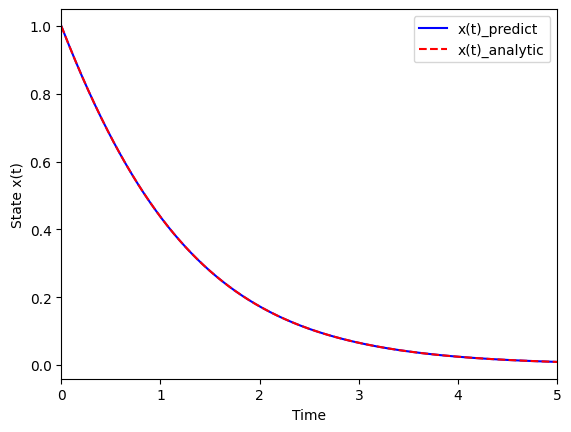}
	\hfil
	\includegraphics[width=0.48\textwidth]{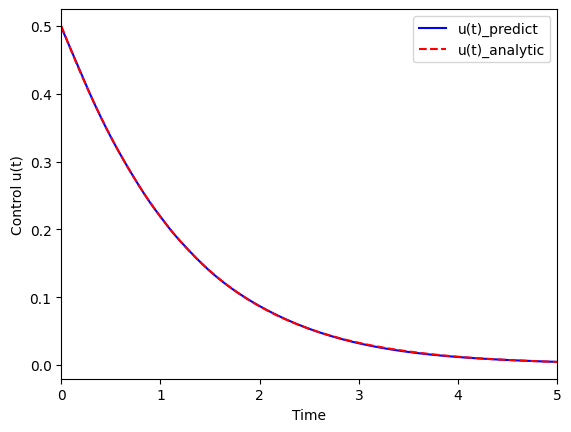}
	\caption{$x(t)$ and $u(t)$ in Example \ref{example of nonlinear} solved by AW-EL-PINNs.}
	\label{fig_exnon:aw_pmpinns}
\end{minipage}
\begin{figure}[H]
\end{figure}
\begin{minipage}{0.48\textwidth}
	\centering
	\includegraphics[width=0.48\textwidth]{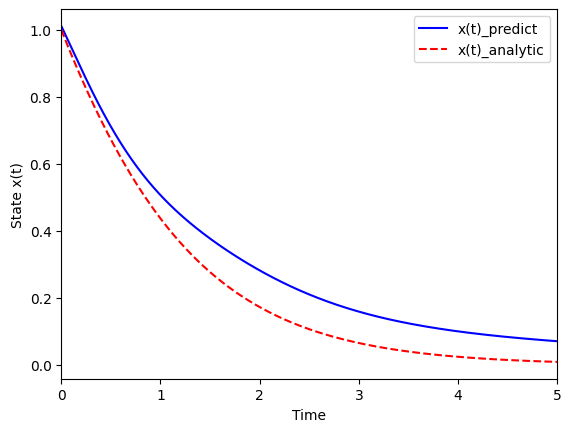}
	\hfil
	\includegraphics[width=0.48\textwidth]{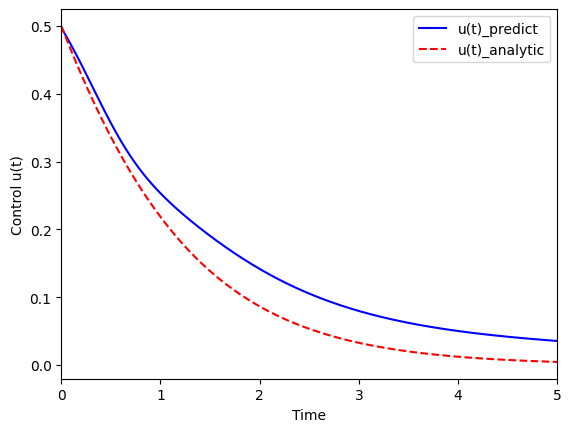}
	\caption{$x(t)$ and $u(t)$ in Example \ref{example of nonlinear} solved by EL-PINNs.}
	\label{fig_exnon:pmpinns}
\end{minipage}
\end{figure}

The update plot of loss weights is shown in Fig. \ref{fig_exnon:exp}. 
\begin{figure}[H]
\centering
\includegraphics[width=0.85\linewidth]{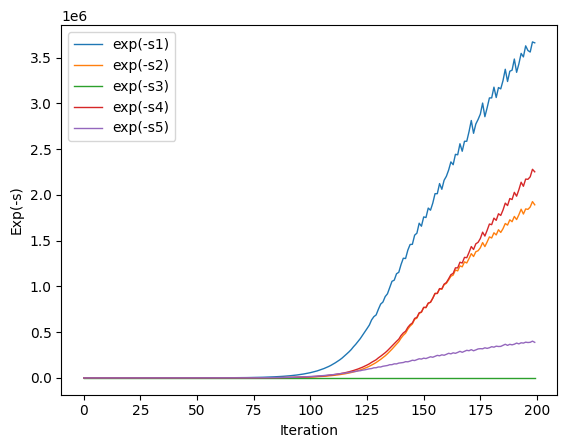}
\caption{$e^{-s_k}$ of AW-EL-PINNs in Example \ref{example of nonlinear}.}
\label{fig_exnon:exp}
\end{figure}

It is evident that the predicted solutions by AW-EL-PINNs closely aligns with the trajectories of analytical solutions. Obviously, the solution predicted by EL-PINNs fails to fit the trajectory of the analytical solutions. 
Hence, it can be conclusively demonstrated that 
AW-EL-PINNs framework exhibits superior predictive accuracy and enhanced stability compared to the conventional EL-PINNs methodology when addressing the problem under consideration.

Notwithstanding the exclusion of  precision comparisons by graphs, we have  conducted comparative convergence analysis experiments.
Fig. \ref{fig_exnon:l2x} and  Fig. \ref{fig_exnon:l2u} illustrate that AW-EL-PINNs converge within approximately 6000 iterations, whereas EL-PINNs fail to converge near 0 even after 20000 iterations, as shown by relative $\mathcal{L}_2$ errors.  
\begin{figure}[H]

\centering
\includegraphics[width=0.85\linewidth]{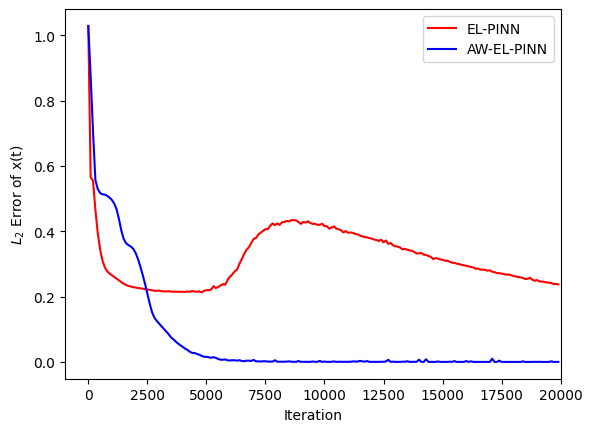}
\caption{Evolution of $x(t)$'s relative $\mathcal{L}_2$ errors in Example \ref{example of nonlinear}.}
\label{fig_exnon:l2x}
\end{figure}
\begin{figure}[H]

\centering
\includegraphics[width=0.85\linewidth]{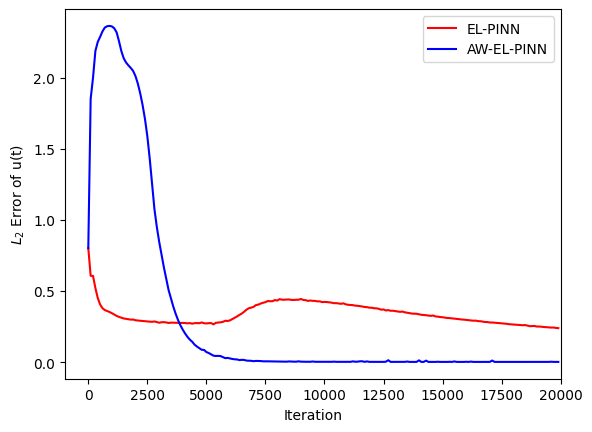}
\caption{Evolution of $u(t)$'s relative $\mathcal{L}_2$ errors in Example \ref{example of nonlinear}.}
\label{fig_exnon:l2u}

\end{figure}

Likewise, we provide the errors obtained from solving this problem using gradient method, conjugate method  and shooting method. The results are illustrated in Table \ref{duibi3}.

\begin{table}[H]  
\centering  
\begin{tabular}{ccc}  
	\hline
	\hline
	Method& MAE & RLE \\
	\hline
	Gradient & 1.25e-02 & 2.55e-02 \\
	Conjugate Gradient & 3.46e-03 & 4.71e-03\\
	Shooting & 5.69e-03 & 4.77e-03  \\
	EL-PINNs & 8.24e-02& 2.39e-01\\
	AW-EL-PINNs & {\bf9.63e-04} &{\bf1.29e-03}\\
	\hline
	\hline
\end{tabular}
\caption{Maximum absolute error and relative $\mathcal{L}_2$ error of numerical methods in Example \ref{example of nonlinear}}  
\label{duibi3}  
\end{table}
AW-EL-PINNs  achieve the highest accuracy in this problem, with both the maximum absolute error  and relative $\mathcal{L}_2$ error
significantly lower than those of other methods. Its superior precision and stability make it the most effective choice. 

\begin{remark}
The observed persistent decline in the relative $\mathcal{L}_2$ error of EL-PINNs demonstrates that with extended training iterations, these errors will asymptotically approach 0. Therefore, 
we increase the training iterations for EL-PINNs to 60000 and  subsequently obtain new trajectories for $x(t)$ and $u(t)$. In addition, the evolution of relative $\mathcal{L}_2$ errors for $x(t)$ and $u(t)$  are also provided.

\begin{figure}[H]
	\centering
	\begin{minipage}{0.48\textwidth}
		\centering
		\includegraphics[width=0.48\textwidth]{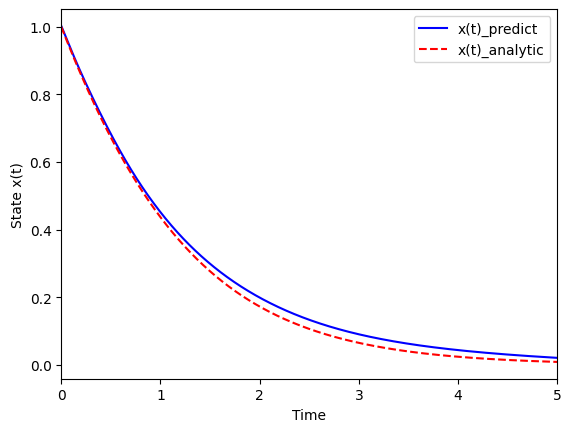}
		\hfil
		\includegraphics[width=0.48\textwidth]{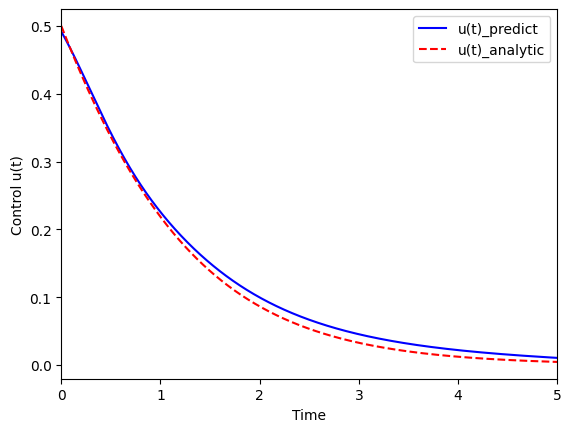}
		\caption{$x(t)$ and $u(t)$ in Example \ref{example of nonlinear} solved by EL-PINNs in 60000 iterations.}
		\label{fig_exnon:pmpinns60000}
	\end{minipage}
\end{figure}
\begin{figure}[H]
	\centering
	\includegraphics[width=0.85\linewidth]{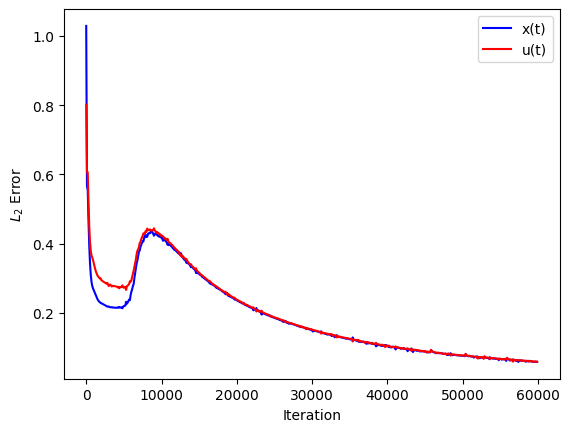}
	\caption{relative $\mathcal{L}_2$ errors in Example \ref{example of nonlinear} solved by EL-PINNs in 60000 iterations.}
	\label{fig_exnash:l2_pmpinns60000}
\end{figure}

As illustrated in Fig. \ref{fig_exnon:pmpinns60000}, compared to Fig. \ref{fig_exnon:aw_pmpinns} (with 20000  iterations), the predictions of EL-PINNs show greater alignment with the analytical solutions and therefore suggest improved accuracy with an increased number of iterations. However, the accuracy of EL-PINNs still falls short when compared to AW-EL-PINNs.

Turning to Fig. \ref{fig_exnash:l2_pmpinns60000}, the trend indicates that as the number of iteration increases, the relative $\mathcal{L}_2$ error moves closer to convergence towards 0, though it has not thoroughly stabilized. This suggests that further extending the training iterations could yield more stable predictions and higher accuracy. Nevertheless, this approach would come at a significant cost in terms of computational resource and time, rendering EL-PINNs less efficient than AW-EL-PINNs.
\end{remark}

\end{example}

\begin{example}\label{2 D nonlinear}(Two-dimensional nonlinear terminal cost  optimal control problem \citep{lasdon2003conjugate})

The Two-dimensional nonlinear optimal control problem is given as follows, for \( t \in [0, 1] \):
\begin{subequations}
\begin{align}
	&\min_uJ(x(t)) = x_2(t_f),  \\
	\text{s.t.}\nonumber\\
	& \dot{x}_1(t) = u(t),  \\
	& x_1(0) = \frac{1}{2}, \\
	& \dot{x}_2(t) = \frac{1}{2}u(t)^2 + u(t)x(t) + u(t) + x_1(t),  \\
	& x_2(0) = 0, 
\end{align}
\end{subequations}
where the analytical solutions are:

\begin{subequations}
\begin{align}
	x_1^*(t) &= \frac{1}{2}t^2 - \frac{3}{2}t + \frac{1}{2},  \\
	x_2^*(t) &= \frac{1}{8}t^4 - \frac{5}{12}t^3 + \frac{3}{8}t^2 - \frac{5}{8}t, \\
	u^*(t) &= t - \frac{3}{2}. 
\end{align}
\end{subequations}

Similar to the previous experiments, we set up four models to solve this problem. For EL-PINNs, the loss function is formulated using 
$\omega_1=\cdots=\omega_8=1$
as the loss weights.

Moreover, in the AW-PINNs model, the loss function is constructed based on $L_J=\frac{1}{N_b}\sum^{N_b}_{j=1}x^j_{2NN}\left(\mathbf{1}[j]\right)$. And  for PINNs, the loss function is defined by incorporating  $\omega_1=\omega_3=\omega_4=\omega_5=1$, $\omega_2=150$ as loss weights.

After training four models for 20000 iterations, the predicted trajectories are plotted below in  Fig. \ref{fig_ex2d:aw_pmpinns}-Fig. \ref{fig_ex2d:pinns}		

\begin{figure}[H]
\centering
\begin{minipage}{0.48\textwidth}
	\centering
	\includegraphics[width=0.48\textwidth]{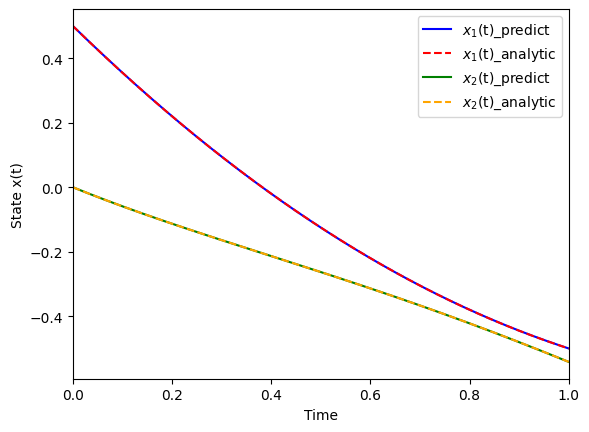}
	\hfil
	\includegraphics[width=0.48\textwidth]{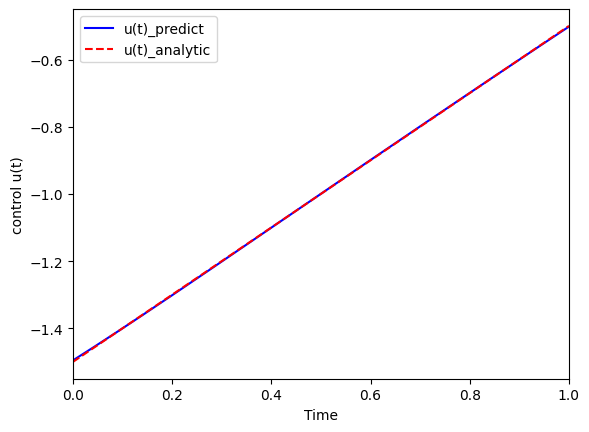}
	\caption{$x(t)$ and $u(t)$ in Example \ref{2 D nonlinear} solved by AW-EL-PINNs.}
	\label{fig_ex2d:aw_pmpinns}
\end{minipage}
\end{figure}
\begin{figure}[H]
\begin{minipage}{0.48\textwidth}
	\centering
	\includegraphics[width=0.48\textwidth]{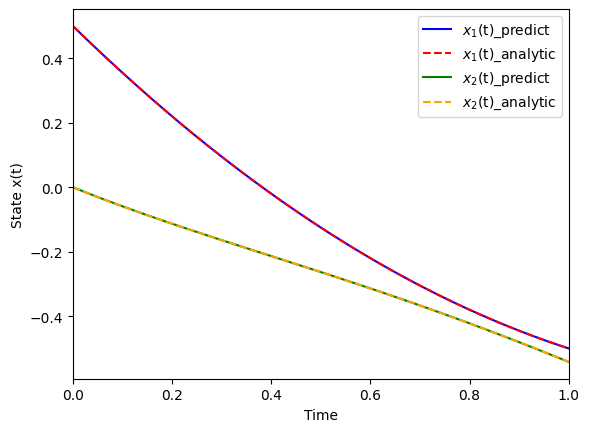}
	\hfil
	\includegraphics[width=0.48\textwidth]{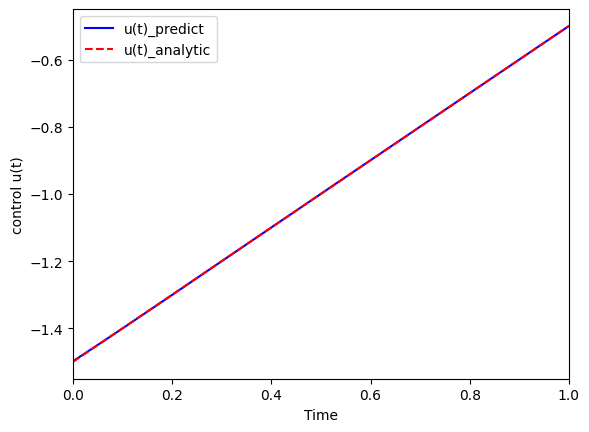}
	\caption{$x(t)$ and $u(t)$ in Example \ref{2 D nonlinear} solved by EL-PINNs.}
	\label{fig_ex2d:pmpinns}
\end{minipage}
\end{figure}
\begin{figure}[H]
\begin{minipage}{0.48\textwidth}
	\centering
	\includegraphics[width=0.48\textwidth]{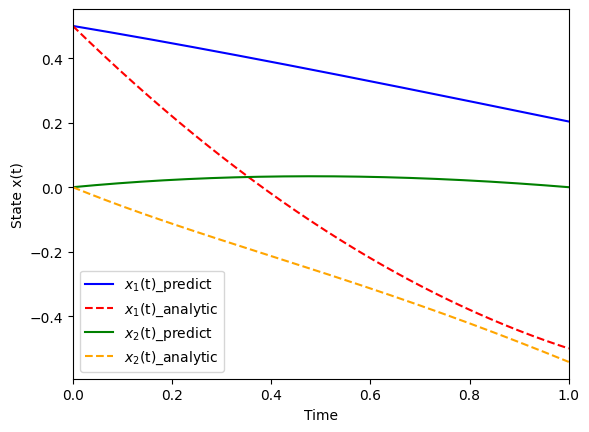}
	\hfil
	\includegraphics[width=0.48\textwidth]{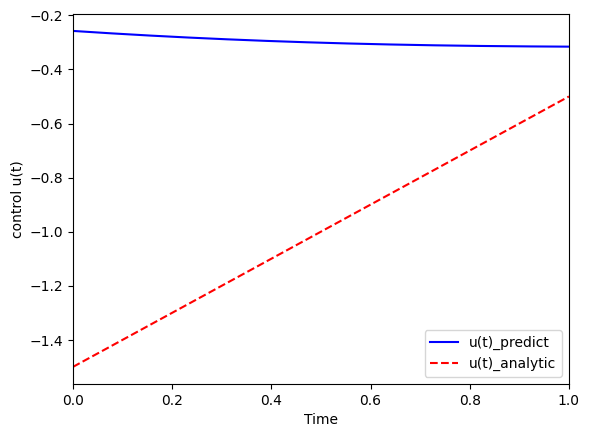}
	\caption{$x(t)$ and $u(t)$ in Example \ref{2 D nonlinear} solved by AW-PINNs.}
	\label{fig_ex2d:aw_pinns}
\end{minipage}
\end{figure}
\begin{figure}[H]
\begin{minipage}{0.48\textwidth}
	\centering
	\includegraphics[width=0.48\textwidth]{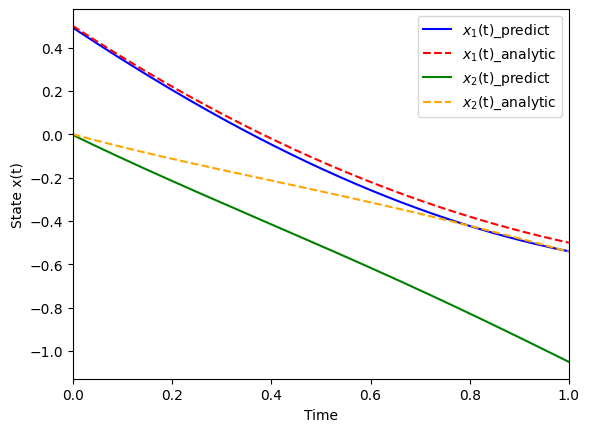}
	\hfil
	\includegraphics[width=0.48\textwidth]{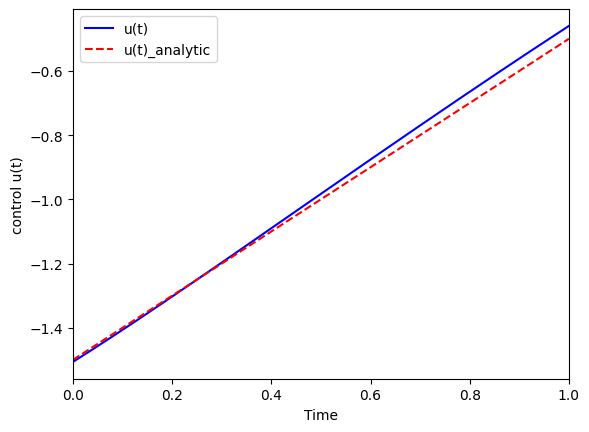}
	\caption{$x(t)$ and $u(t)$ in Example \ref{2 D nonlinear} solved by PINNs.}
	\label{fig_ex2d:pinns}
\end{minipage}
\end{figure}

Fig. \ref{fig_ex2d:aw_pmpinns} and Fig. \ref{fig_ex2d:pmpinns} 
demonstrate that the solutions predicted by AW-EL-PINNs and EL-PINNs align well with the analytical solutions. However, shown in Fig. \ref{fig_ex2d:aw_pinns} and Fig. \ref{fig_ex2d:pinns}, AW-PINNs fail to produce accurate predictions, while PINNs manage to predict $x_1(t)$ and $u(t)$
reasonably well but exhibit poor performance in predicting the solution for $x_2(t)$.

To further compare the prediction accuracy of  AW-EL-PINNs and EL-PINNs, the boxplots of absolute errors are presented as follows:

\begin{figure}[H]
\centering
\includegraphics[width=0.85\linewidth]{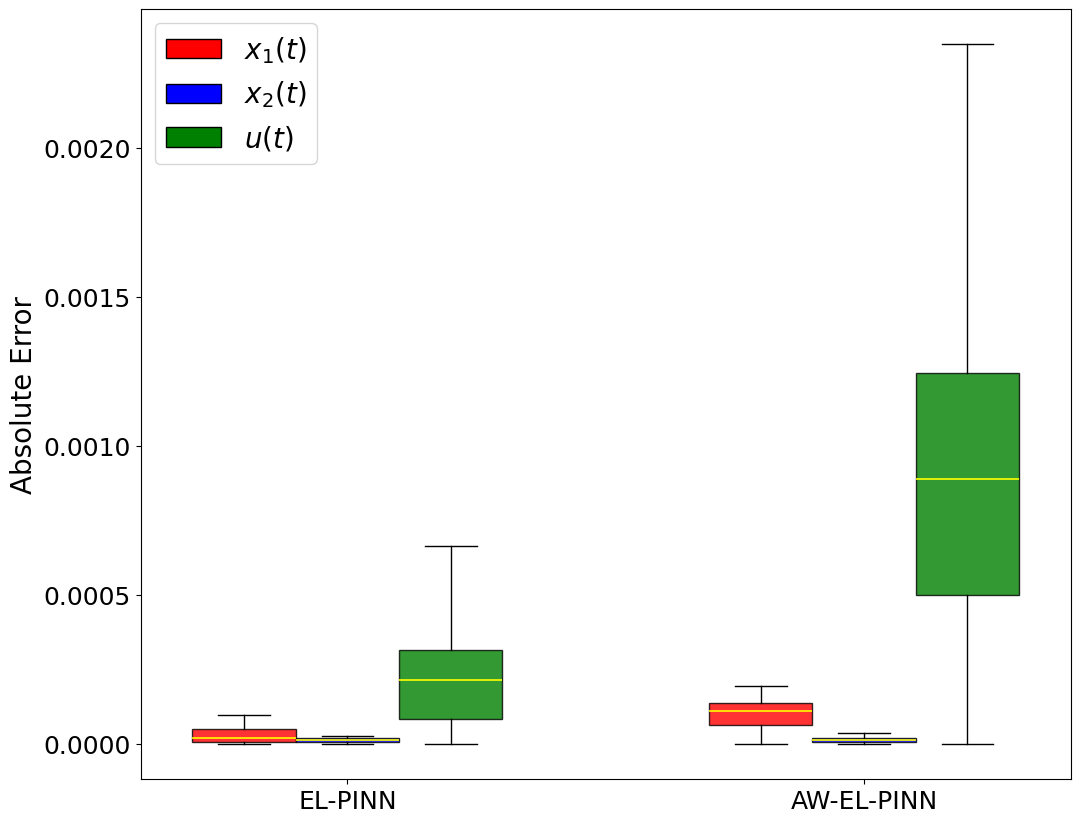}
\caption{Box plots of absolute errors in Example \ref{2 D nonlinear}.}
\label{fig_ex2d:box}
\end{figure}

From the observation of Fig. \ref{fig_ex2d:box}, it can be concluded that EL-PINNs exhibit slightly better prediction accuracy for the three variables compared to AW-PINNs.

In the subsequent step, we also present the bar charts of relative $\mathcal{L}_2$ errors to compare the stability among the models.

\begin{figure}[H]
\centering
\includegraphics[width=0.85\linewidth]{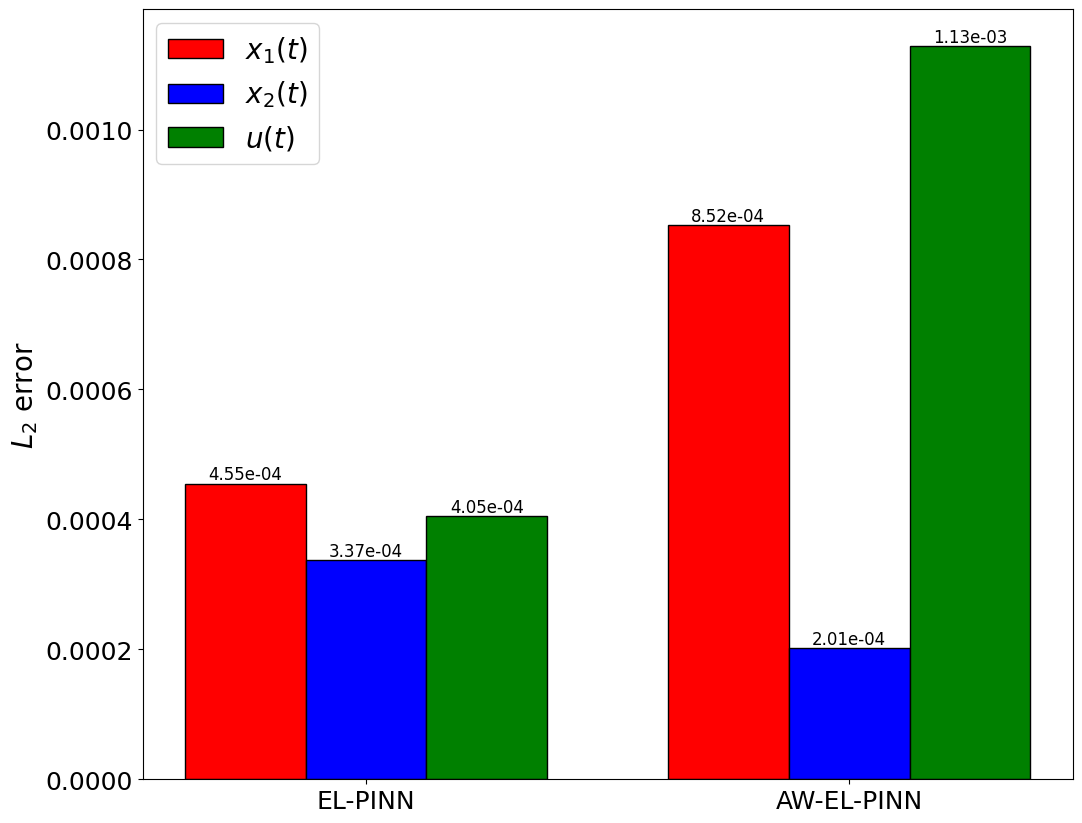}
\caption{Bar charts of relative $\mathcal{L}_2$ errors in Example \ref{2 D nonlinear}.}
\label{fig_ex2d:bar}
\end{figure}

It is evident that, except for 
$x_2(t)$, the relative $\mathcal{L}_2$ errors of 
$x_1(t)$ and $u(t)$ predicted by EL-PINNs are lower than those of AW-EL-PINNs. This indicates that EL-PINNs demonstrate marginally better stability compared to AW-EL-PINNs.

Finally, we need to compare the convergence speeds of the models. The following figures illustrate the relative $\mathcal{L}_2$ errors evolution  over iterations.

\begin{figure}[H]
\centering
\includegraphics[width=0.85\linewidth]{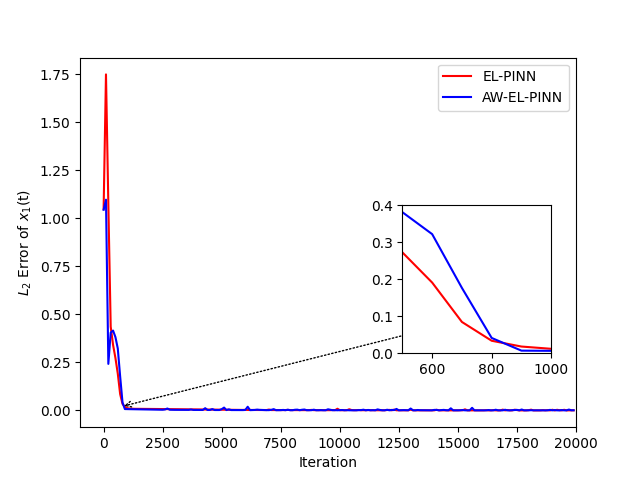}
\caption{Evolution of $x(t)$'s relative $\mathcal{L}_2$ errors in Example \ref{2 D nonlinear}.}
\label{fig_ex2d:xL2}
\end{figure}

\begin{figure}[H]
\centering
\includegraphics[width=0.85\linewidth]{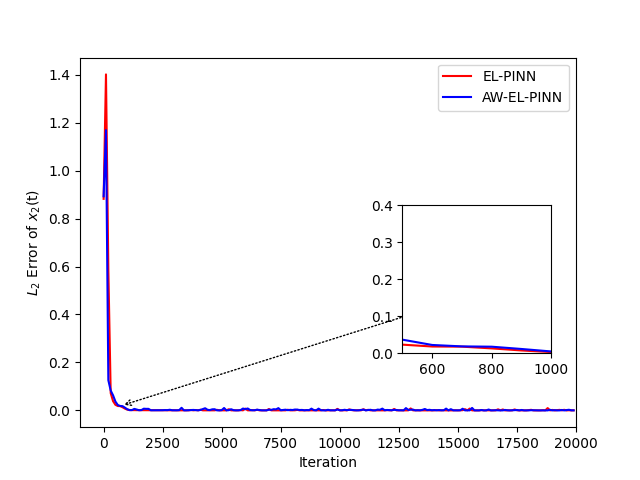}
\caption{Evolution of $u_1(t)$'s relative $\mathcal{L}_2$ errors in Example \ref{2 D nonlinear}.}
\label{fig_ex2d:u1L2}
\end{figure}

\begin{figure}[H]
\centering
\includegraphics[width=0.85\linewidth]{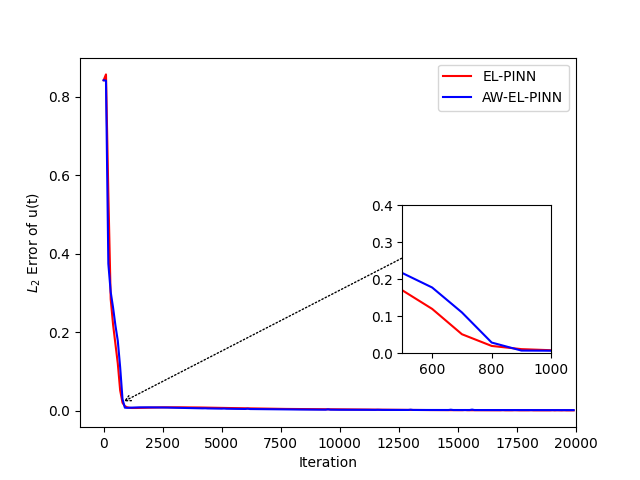}
\caption{Evolution of $u_2(t)$'s relative $\mathcal{L}_2$ errors in Example \ref{2 D nonlinear}.}
\label{fig_ex2d:u2L2}
\end{figure}
As illustrated in Fig. \ref{fig_ex2d:xL2}-Fig. \ref{fig_ex2d:u2L2}, the convergence speeds of AW-EL-PINNs and EL-PINNs are comparable.

In the same manner, we provide the errors obtained from solving this problem using other three conventional  methods. The results are shown in Table \ref{duibi4}.
\begin{table}[h]  
\centering  
\begin{tabular}{ccc}  
	\hline
	\hline
	Method & MAE & RLE   \\
	\hline
	Gradient & 1.90e-01 & 2.41e-01  \\
	Conjugate Gradient &7.11e-02 &1.11e-01\\
	Shooting & 2.57e-01 & 6.74e-02  \\
	EL-PINNs & {\bf 6.93e-04} & {\bf 3.99e-04}\\
	AW-EL-PINNs & 1.56e-03 & 7.27e-04\\
	\hline\hline
\end{tabular}
\caption{Maximum absolute error and relative $\mathcal{L}_2$ error of numerical methods in Example \ref{2 D nonlinear}}  
\label{duibi4}  
\end{table}

AW-EL-PINNs and EL-PINNs  achieve the higher accuracy, with both kinds of errors significantly lower than those of gradient and shooting, where EL-PINNs shows more accurate that AW-EL-PINNs.

\begin{remark}
From Fig. \ref{fig_ex2d:aw_pinns}, it can be observed that AW-PINNs fail to effectively solve this problem. We experimented with various methods to improve its performance and found that fixing the loss weight of  $L_5$ as $\omega_5=1$, while keeping adaptive weights for the other 
$L_k(k=1,2,3,4)$, yields satisfactory prediction results, as shown by Fig. \ref{fig_ex2d:awpinns_w5=1} and the relative $\mathcal{L}_2$ errors evolution by Fig. \ref{fig_ex2d:awpinns_w5=1_L2}.

\begin{figure}[H]
	\centering
	\begin{minipage}{0.48\textwidth}
		\centering
		\includegraphics[width=0.48\textwidth]{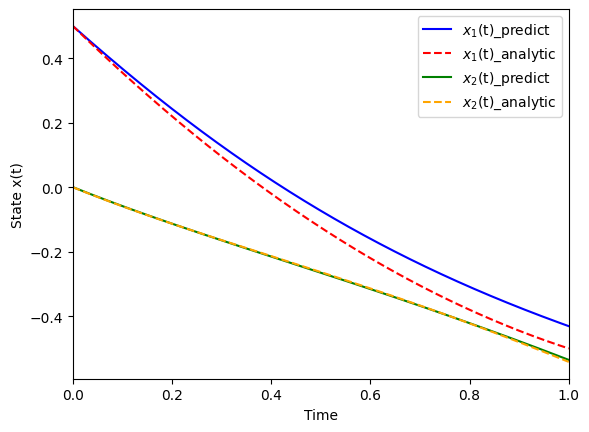}
		\hfil
		\includegraphics[width=0.48\textwidth]{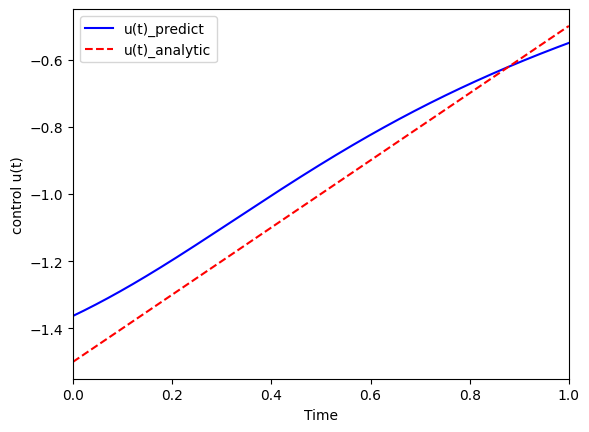}
		\caption{$x(t)$ and $u(t)$ in Example \ref{2 D nonlinear} solved by AW-PINNs with $\omega_5=1$.}
		\label{fig_ex2d:awpinns_w5=1}
	\end{minipage}
\end{figure}

\begin{figure}[H]
	\centering
	\includegraphics[width=0.85\linewidth]{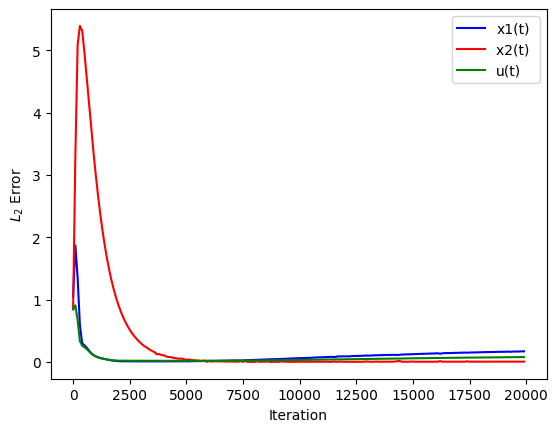}
	\caption{relative $\mathcal{L}_2$ errors in Example \ref{2 D nonlinear} solved by AW-PINNs with $\omega_5=1$.}
	\label{fig_ex2d:awpinns_w5=1_L2}
\end{figure}

It is evident that the prediction accuracy has improved significantly compared to that shown in Fig. \ref{fig_ex2d:aw_pinns}. However, it still falls far short of the prediction accuracy achieved by AW-EL-PINNs and EL-PINNs.

\end{remark}
\end{example}

\begin{example}\label{3o}(Nonlinear damped oscillator energy consumption and terminal state control.)

In this example, we consider the following idealized cubic damping system:

\begin{subequations}
\begin{align}
	& \dot{x}(t)=v(t),  \\
	& x(0)=1 \textnormal{m} \\
	& \dot{v}(t)=-\frac{k}{m}x(t)-\frac{c}{m}v(t)^3+\frac{u(t)}{m},  \\
	& v(0) = -0.8\textnormal{m/s}, 
\end{align}
\end{subequations}
where  $x(t)$ represents the displacement of the oscillator, and $v(t)$ represents its velocity. The linear feedback control gain $k$ is taken as $2\,\text{N/m}$, and the cubic damping coefficient is $c = 0.005\,\text{N}\cdot\text{s}^3/\text{m}^3$. The external force acting on the system is denoted by $u(t)$, with units of $\text{N}$. The mass $m$ is taken as $1\,\text{kg}$.
The value function is defined as follows:
\begin{align}V(u(t),x(t),v(t))&=\min_u J(u(t),x(t),v(t))\notag\\&=\min_u\int^{3}_0\frac{1}{2}u(t)^2\mathrm{d}t+25x(3)^2+25v(3)^2,\end{align}
which implies minimizing the energy consumption within 3$s$ while controlling the terminal state to the origin.

Since this problem is a nonlinear optimal control problem with no analytical solutions, we employ the fourth-order Runge-Kutta (RK-4) method which solves the differential equations of TPBVP  for this problem  as the  reference to approximate the exact solutions. Next, we set the learning rates of both optimizers in AW-EL-PINNs to $1e-4$. All loss weights in EL-PINNs are set to 1. Both networks are trained for 100000 iterations, and the resulting outcomes are illustrated in Fig .\ref{fig:3oaw}-Fig. \ref{fig:3ohand}.

\begin{figure}[H]
\begin{minipage}{0.48\textwidth}
	\centering
	\includegraphics[width=0.48\textwidth]{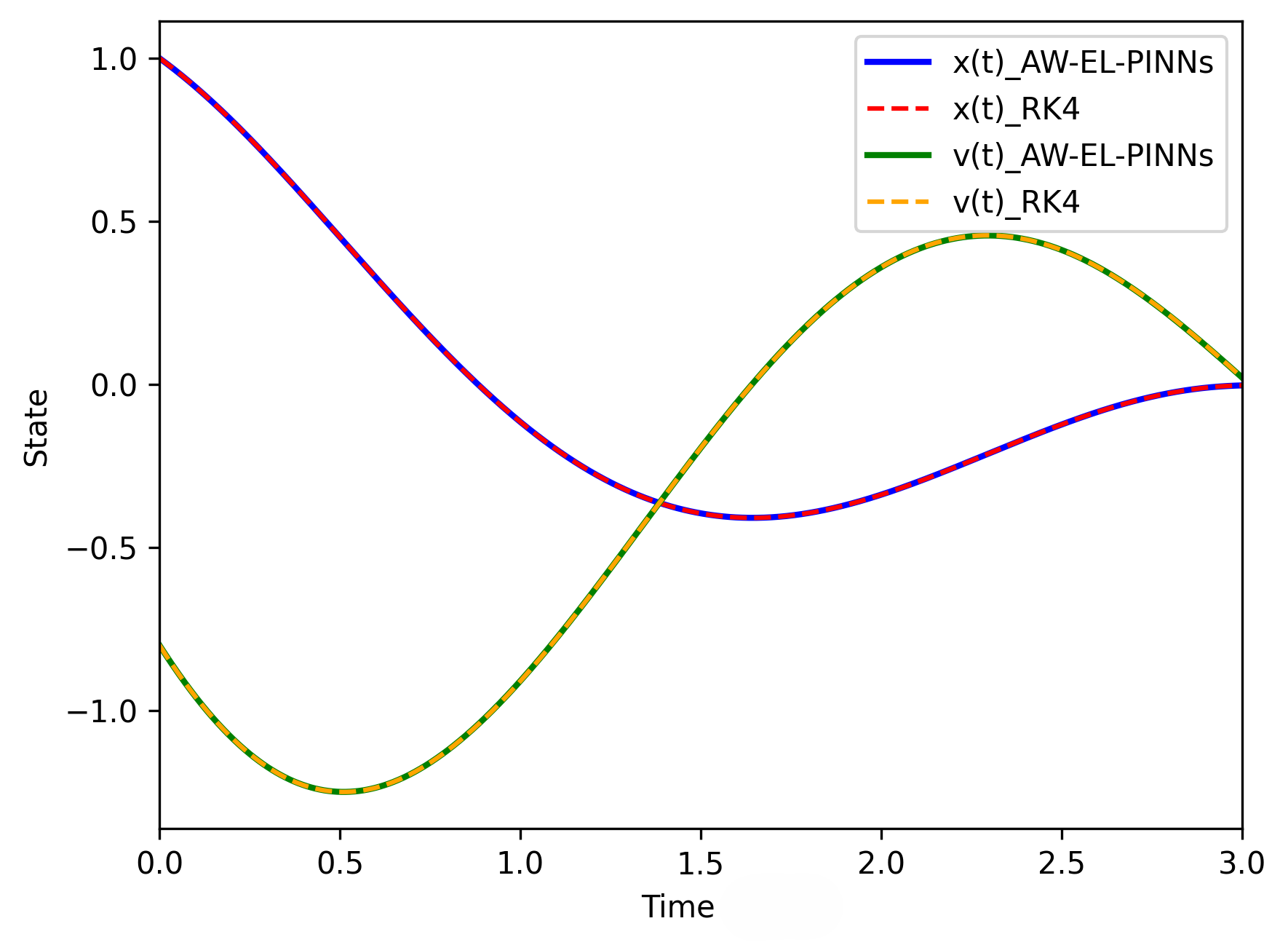}
	\hfil
	\includegraphics[width=0.48\textwidth]{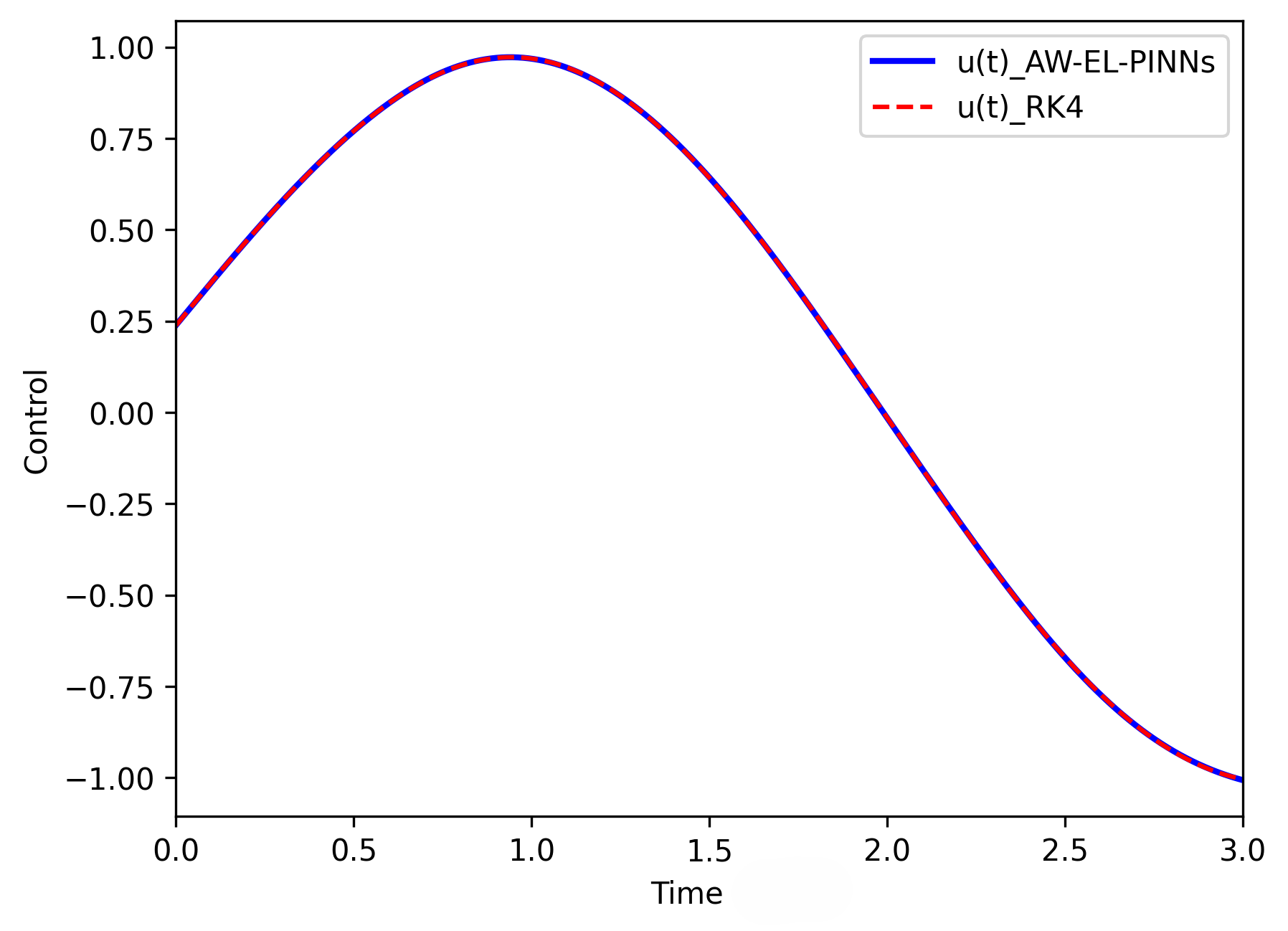}
	\caption{$x(t)$ and $u(t)$ in Example \ref{3o} solved by AW-EL-PINNs.}
	\label{fig:3oaw}
\end{minipage}
\end{figure}
\begin{figure}[H]
\begin{minipage}{0.48\textwidth}
	\centering
	\includegraphics[width=0.48\textwidth]{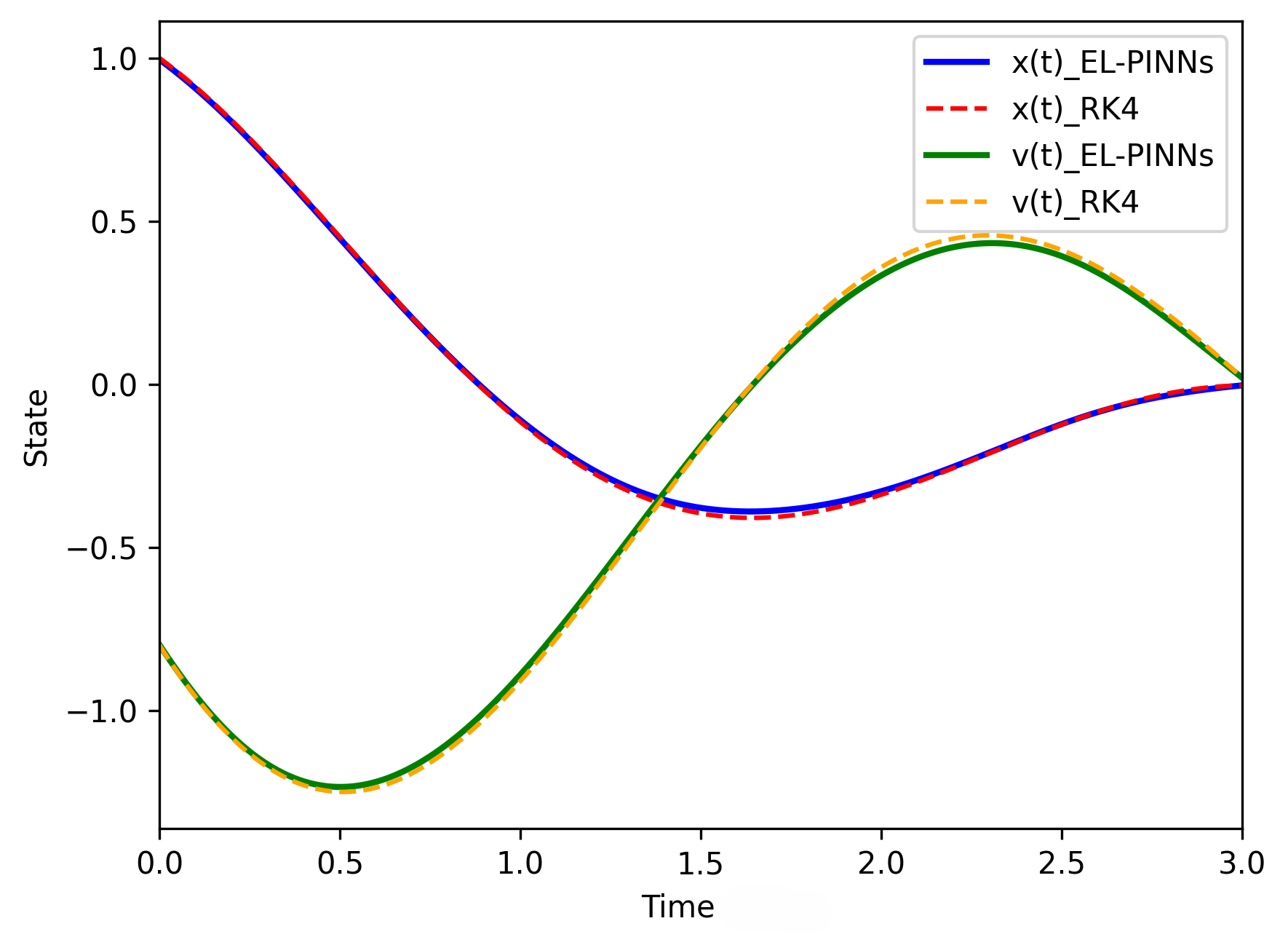}
	\hfil
	\includegraphics[width=0.48\textwidth]{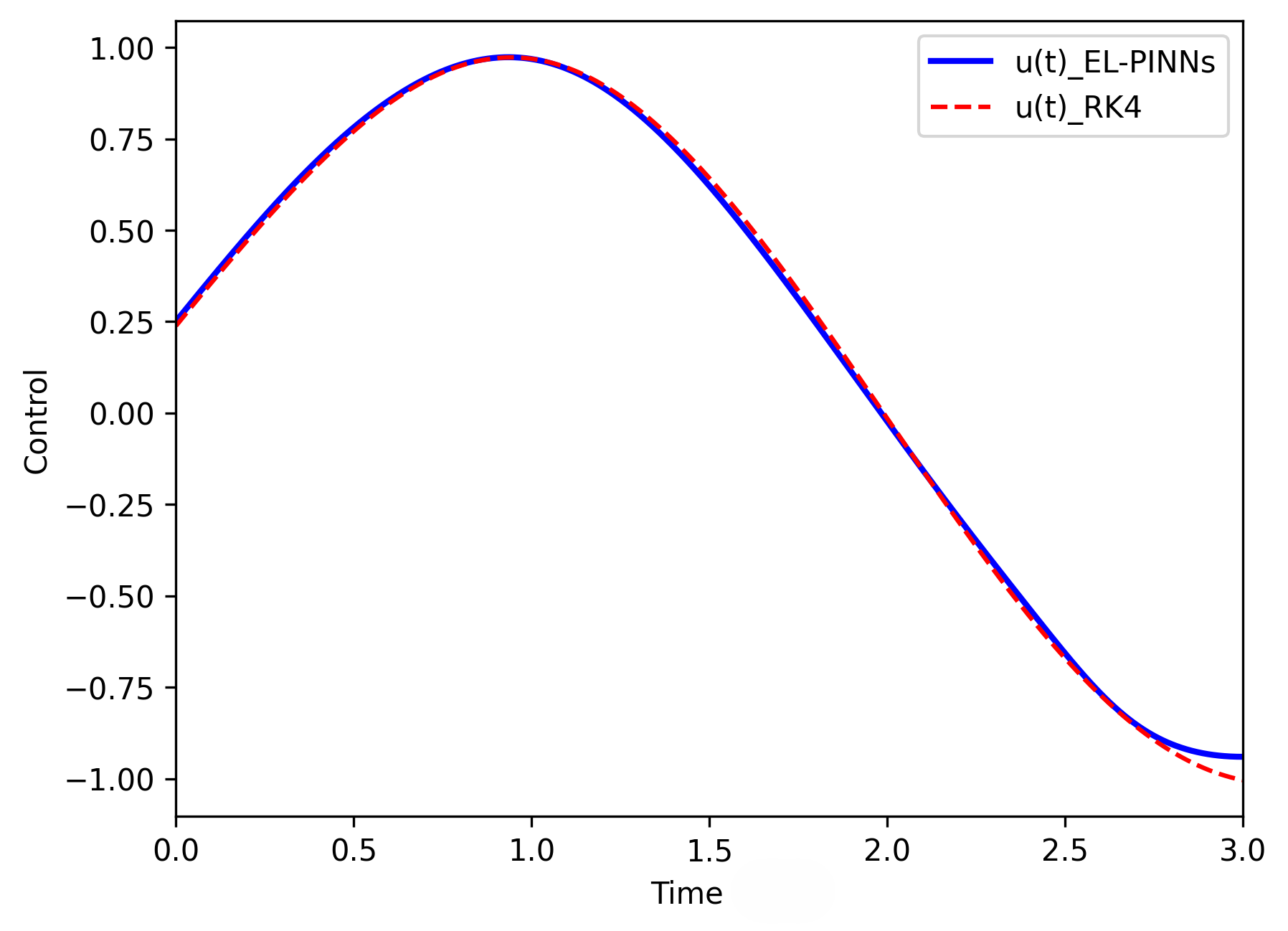}
	\caption{$x(t)$ and $u(t)$ in Example \ref{3o} solved by EL-PINNs.}
	\label{fig:3ohand}
\end{minipage}
\end{figure}

It is evident that both AW-EL-PINNs and EL-PINNs are capable of accurately predicting the trajectories of state variables and control variable. Moreover, it is obvious that the prediction accuracy of trajectories using AW-EL-PINNs is superior to that of EL-PINNs. Therefore, we do not separately plot box plots and bar charts to discuss the comparison of accuracy and stability between the two algorithms. Instead, we present specific numerical values in the subsequent list-based comparison.

Next, we compare the gradient method, conjugate gradient method, AW-EL-PINNs, and EL-PINNs through table \ref{duibi5}. The shooting method is excluded from this comparison due to significant deviations in its computational results. 	\begin{table}[h]  
\centering  
\begin{tabular}{ccc}  
	\hline
	\hline
	Method & MAE & RLE   \\
	\hline
	Gradient & 2.06e-2  &  1.59e-2  \\
	Conjugate Gradient &{\bf 9.67e-4 }  & 9.73e-4\\
	EL-PINNs &3.69e-2   & 2.37e-2\\
	AW-EL-PINNs &1.20e-3 &   {\bf 5.86e-4}\\
	\hline\hline
\end{tabular}
\caption{Maximum absolute error and relative $\mathcal{L}_2$ error of numerical methods in Example \ref{3o}}  
\label{duibi5}  
\end{table}

Consequently, AW-EL-PINNs achieves the optimal relative $ \mathcal{L}_2$ error while maintaining the second-lowest maximum absolute error, demonstrating 1-2 orders of magnitude improvement over traditional gradient methods and significant superiority over EL-PINNs. 

\end{example}

\begin{example}\label{car}
(Energy consumption and terminal state control of cart motion under nonlinear horizontal friction)

Assume a cart undergoes decelerated horizontal motion on a rough surface under an applied thrust. The smoothness of the surface varies with displacement, becoming smoother as the displacement increases in the positive direction. Let the positive direction be to the right. The cart starts its decelerated motion in the negative direction from a position of 2 m, with an initial velocity of -0.8 m/s. The dynamic system governing the motion is given as follows:
\begin{subequations}
	\begin{align}
		& \dot{x}(t)=v(t),  \\
		& x(0)=2 \textnormal{m} \\
		& \dot{v}(t)=u(t)-c_0\exp(-kx(t))v(t)\sqrt{v(t)^2+\varepsilon^2}, \label{35c} \\
		& v(0) = -0.8\textnormal{m/s}. 
	\end{align}
\end{subequations}
\end{example}
Here, $c_0 = 0.01 \, \text{m}^{-1}$ represents the roughness coefficient at $x = 2$ m, while $k = -0.005 \, \text{m}^{-1}$ denotes the exponential decay coefficient of the surface roughness. $u(t)$ represents the acceleration generated by the applied thrust, whereas $x(t)$ and $v(t)$ denote the displacement and velocity, respectively.  
\eqref{35c} describes the time derivative of velocity, which equals the acceleration due to thrust minus the acceleration induced by friction. Since the friction force always acts in the direction opposite to velocity, \eqref{35c} should ideally be written as
 $\dot{v}(t)=u(t)-c_0\exp(-kx(t))v(t)|v(t)|$.  
However, considering that $|v(t)|$ is not differentiable everywhere with respect to $v(t)$, a regularization term $\varepsilon$ is introduced, where $\varepsilon$ is set to $1e-6$ $ m/s$.

The value function is defined as follows:
\begin{align}V(u(t),x(t),v(t))&=\min_u J(u(t),x(t),v(t))\notag\\&=\min_u\int^{5}_0\frac{1}{2}u(t)^2\mathrm{d}t+25x(5)^2+25v(5)^2,\end{align}
which indicates minimizing the energy consumption within $5s$ while controlling the terminal state to the origin.

\begin{figure}[H]
	\centering
	\includegraphics[width=0.85\linewidth]{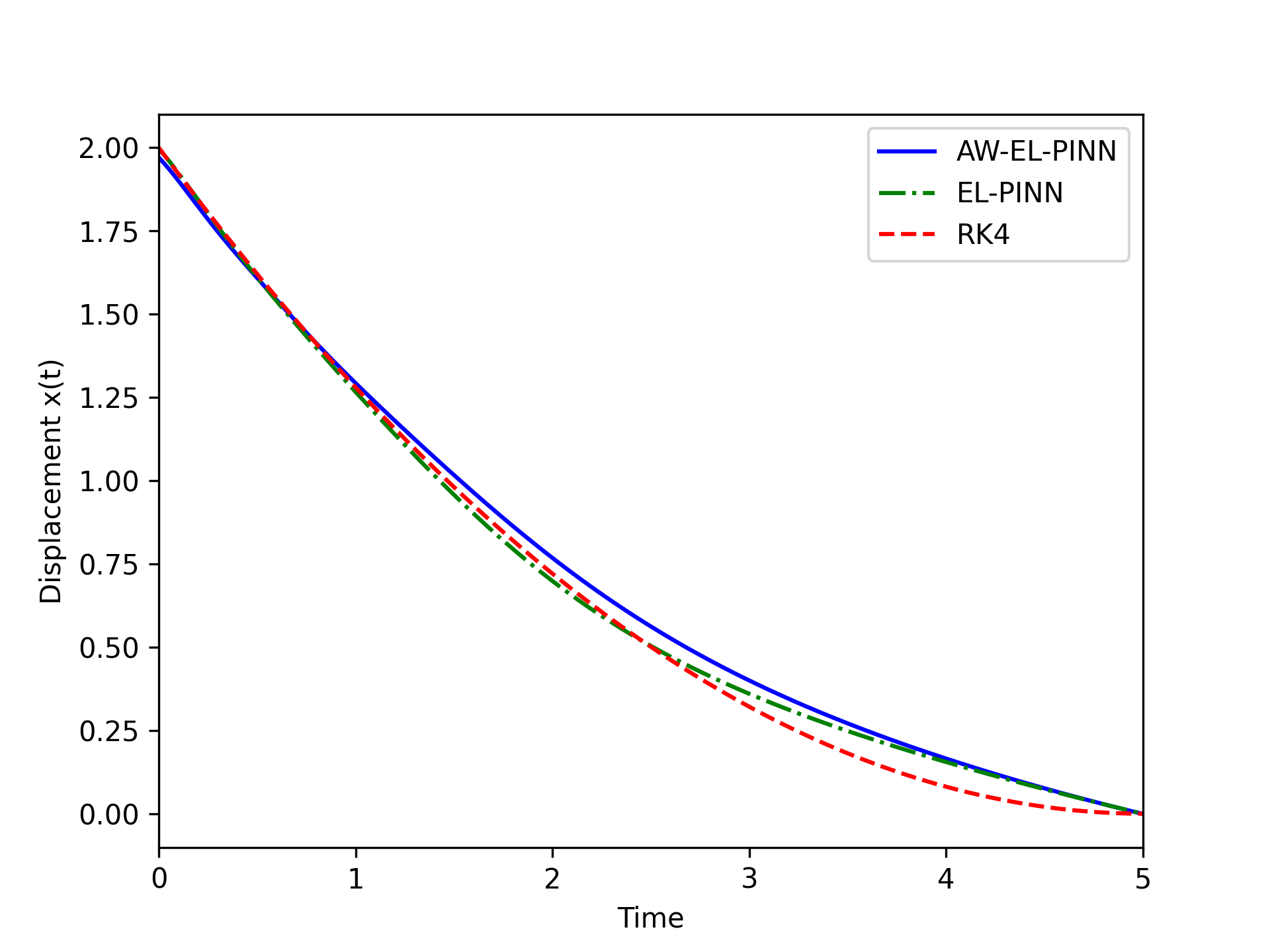}
	\caption{The  trajectories of displacement $x(t)$ in Example \ref{car}.}
	\label{fig_5x}
\end{figure}

\begin{figure}[H]
	\centering
	\includegraphics[width=0.85\linewidth]{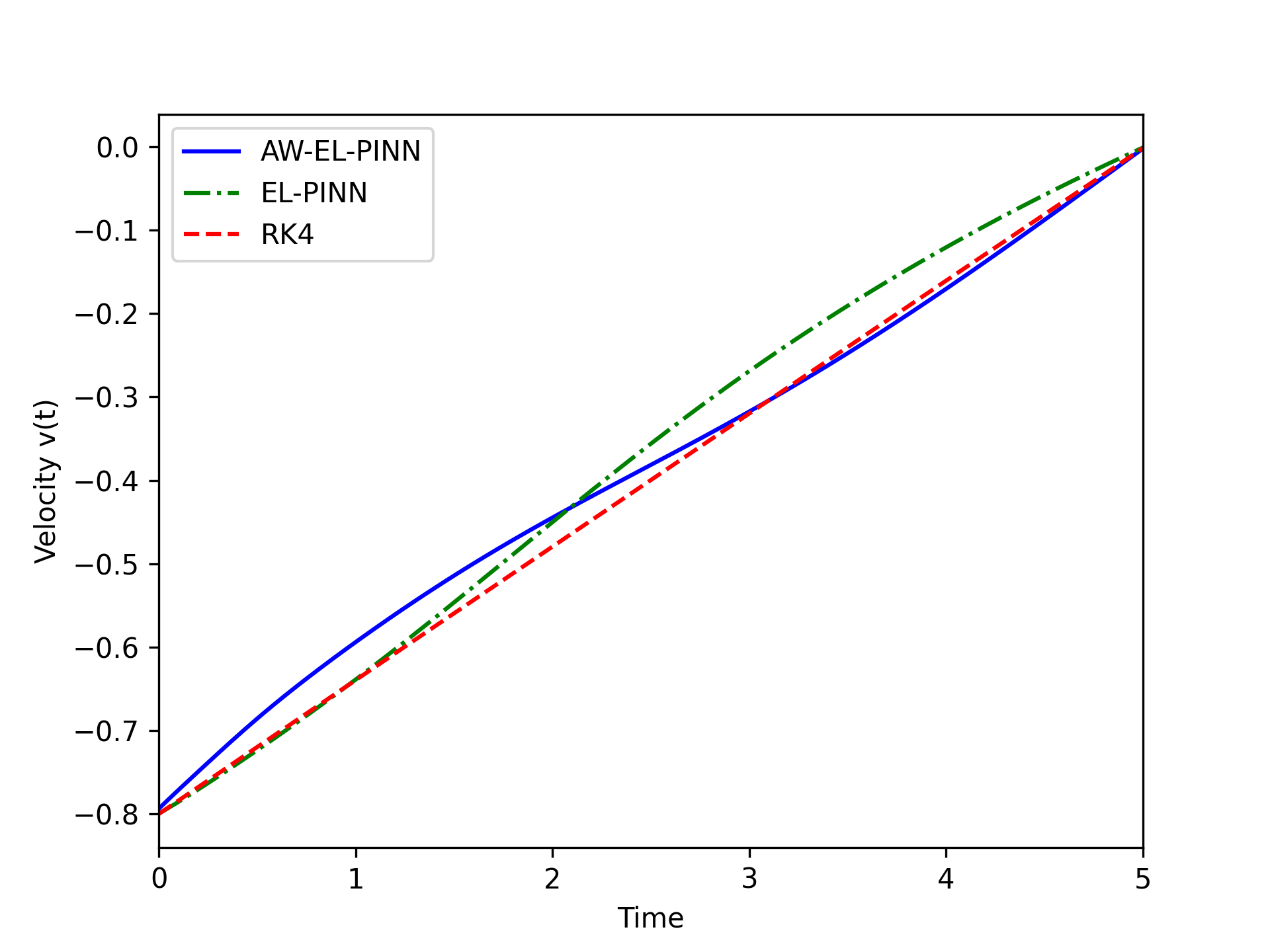}
	\caption{The  trajectories of velocity $v(t)$ in Example \ref{car}.}
	\label{fig_5v}
\end{figure}

\begin{figure}[H]
	\centering
	\includegraphics[width=0.85\linewidth]{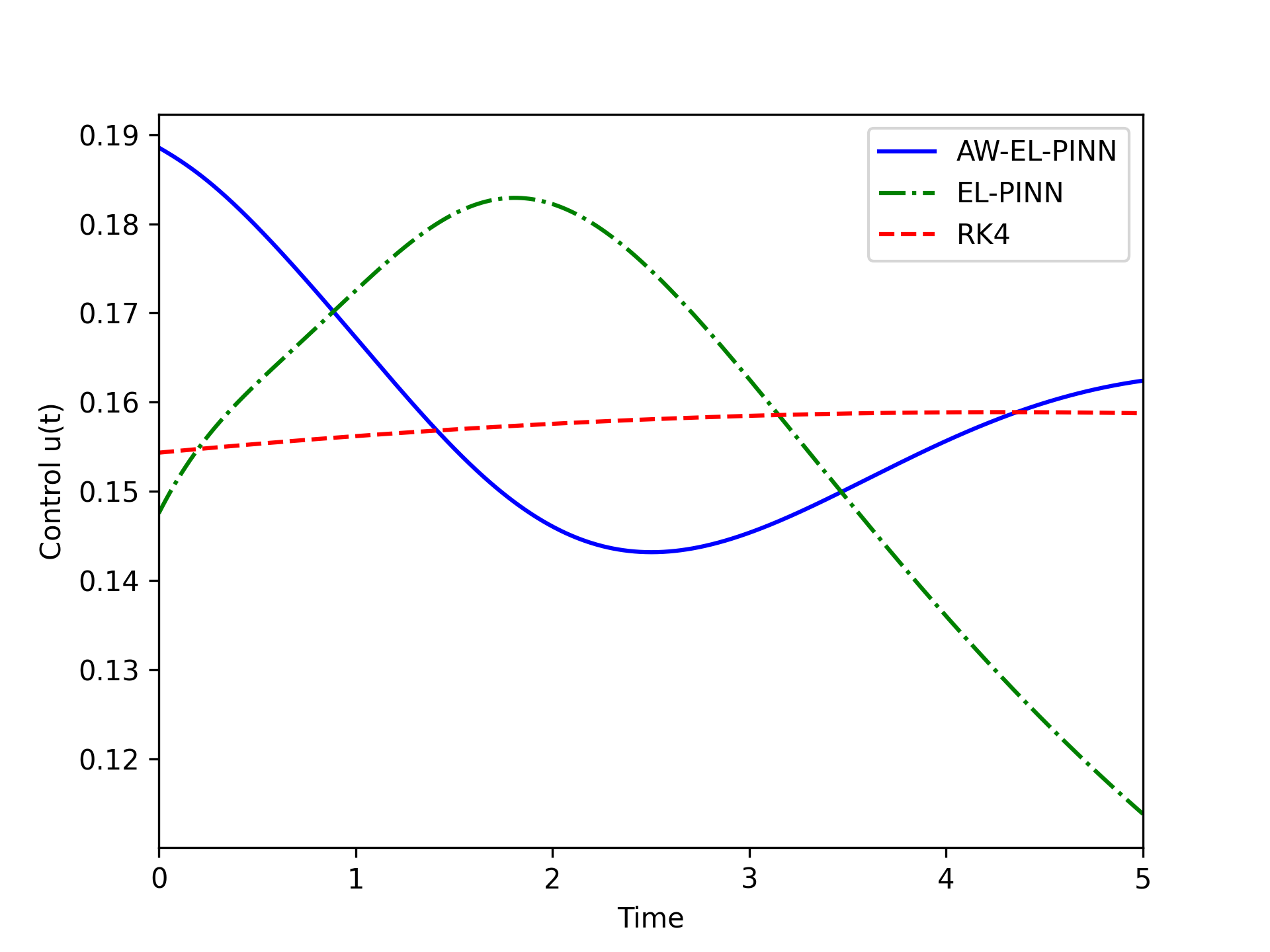}
	\caption{The  trajectories of control $u(t)$ in Example \ref{car}.}
	\label{fig_5u}
\end{figure}

The analytical solutions do not exist for this example either. Therefore, the RK-4 algorithm is used to solve the TPBVP, which represents for exact solutions. In this example, only AW-EL-PINNs and EL-PINNs are employed for the solution, with the number of iterations set to 30000. The loss weights for EL-PINNs are all set to 1. The resulting trajectories of the variables are shown in Fig. \ref{fig_5x}-Fig. \ref{fig_5u}.

From the figures, it can be observed that both AW-EL-PINNs and EL-PINNs can accurately predict the displacement and velocity of the cart. However, the predicted trajectories of the control variables differ significantly. In the RK-4 reference solution, $u(t)$ is explicitly reconstructed from the adjoint variable $\lambda(t)$, whereas in AW-EL-PINNs the control
	$u(t)$ appears only indirectly through the Hamiltonian. This lack of direct supervision leads to the pronounced discrepancy between the predicted control trajectory and the reference in this example.

\section{Conclusion}\label{Sec:conclusion}

This study introduces AW-EL-PINNs, which integrate the Euler-Lagrange theorem with PINNs to establish a framework that converts  Euler-Lagrange systems in optimal control problems into numerically solvable ODEs, which provides a methodological approach for strategy acquisition in potential physical or engineering applications.

Comparative experiments yield several important findings. Most notably, AW-EL-PINNs consistently exhibit superior predictive performance compared to  PINNs in \citet{mowlavi2023optimal}, due to their adaptive weighting capability and avoiding the complexities associated with  selection of sub-loss function with $J$. Additionally, while optimal loss weights could theoretically enable EL-PINNs to match AW-EL-PINNs' performance, finding these weights significantly increases experimental complexity, reinforcing AW-EL-PINNs as a practical and efficient alternative. Furthermore, AW-EL-PINNs demonstrate remarkable accuracy, stability. However, the accuracy achieved in predicting the control variable $u(t)$ is generally lower than that for the state variable $x(t)$, primarily because $x(t)$ benefits from supervision via initial conditions, whereas $u(t)$  is solved exclusively from residual equations.

All examples in this paper focus on unconstrained control problems. For the common box constraint $u_{\min} \leq u(t) \leq u_{\max}$, feasibility can be enforced structurally by applying a bounded activation function to the network output, such as tanh or sigmoid, which ensures that the predicted control  stays within the prescribed limits. We intend to investigate these strategies  in future work.

Looking ahead, we aim to extend the framework of PINNs to the solutions of differential games which focus on decision making processes involving multiple agents, with the goal of achieving coordinated optimization of control variables presented by each agent while balancing individual and collective objectives. Although existing studies on evolutionary games  primarily focus on discrete strategy spaces \citep{wang2015universal,wang2015evolutionary,wang2017onymity,wang2018exploiting,li2018punishment,
	wang2022modelling}, they share deep commonalities with continuous-state differential games—both require a trade-off between individual rationality and collective welfare. Based on this insight, we are going to propose the multi-task learning framework that incorporates physical constraints and adaptive weight allocation, offering a novel paradigm for obtaining numerical solutions of differential games.

On the other hand, we will focus on developing numerical solution algorithms for impulsive optimal  control and impulsive differential games. Impulsive  control, characterized by abrupt changes in state variables, presents significant challenges due to the need for precise handling of impulsive dynamical systems \citep{Qiu2022Fuzzy,
	Wang2024State,
	Wang2024EventTriggered,
	wang2024Stability,
	he2025Synchronization,
	wang2025synchronisation}. Besides, impulsive differential games concentrate on optimizing the timing and magnitude of impulses while addressing the computation of equilibrium solutions.
Building upon existing theoretical foundations \citep{chahim2012tutorial,tauchnitz2015pontryagin,sadana2022feedback,lv2022nonzero,Li2024Verification} and integrating  with the methodologies of PINNs, we aim to devise advanced numerical algorithms tailored to these problems.

\bibliography{gobib}

\begin{thebibliography}{50}
\expandafter\ifx\csname natexlab\endcsname\relax\def\natexlab#1{#1}\fi
\providecommand{\url}[1]{\texttt{#1}}
\providecommand{\href}[2]{#2}
\providecommand{\path}[1]{#1}
\providecommand{\DOIprefix}{doi:}
\providecommand{\ArXivprefix}{arXiv:}
\providecommand{\URLprefix}{URL: }
\providecommand{\Pubmedprefix}{pmid:}
\providecommand{\doi}[1]{\href{http://dx.doi.org/#1}{\path{#1}}}
\providecommand{\Pubmed}[1]{\href{pmid:#1}{\path{#1}}}
\providecommand{\bibinfo}[2]{#2}
\ifx\xfnm\relax \def\xfnm[#1]{\unskip,\space#1}\fi
\bibitem[{Antonelo et~al.(2024)Antonelo, Camponogara, Seman, Jordanou, de~Souza
  \& H{\"u}bner}]{antonelo2024physics}
\bibinfo{author}{Antonelo, E.~A.}, \bibinfo{author}{Camponogara, E.},
  \bibinfo{author}{Seman, L.~O.}, \bibinfo{author}{Jordanou, J.~P.},
  \bibinfo{author}{de~Souza, E.~R.}, \& \bibinfo{author}{H{\"u}bner, J.~F.}
  (\bibinfo{year}{2024}).
\newblock \bibinfo{title}{Physics-informed neural nets for control of dynamical
  systems}.
\newblock {\it \bibinfo{journal}{Neurocomputing}\/},  {\it
  \bibinfo{volume}{579}\/}, \bibinfo{pages}{127419}.
\bibitem[{Chahim et~al.(2012)Chahim, Hartl \& Kort}]{chahim2012tutorial}
\bibinfo{author}{Chahim, M.}, \bibinfo{author}{Hartl, R.~F.}, \&
  \bibinfo{author}{Kort, P.~M.} (\bibinfo{year}{2012}).
\newblock \bibinfo{title}{A tutorial on the deterministic impulse control
  maximum principle: necessary and sufficient optimality conditions}.
\newblock {\it \bibinfo{journal}{European Journal of Operational Research}\/},
  {\it \bibinfo{volume}{219}\/}, \bibinfo{pages}{18--26}.
\bibitem[{Dong et~al.(2025)Dong, Zhu, An, Jiang \& Ma}]{dong2025barrier}
\bibinfo{author}{Dong, B.}, \bibinfo{author}{Zhu, X.}, \bibinfo{author}{An,
  T.}, \bibinfo{author}{Jiang, H.}, \& \bibinfo{author}{Ma, B.}
  (\bibinfo{year}{2025}).
\newblock \bibinfo{title}{Barrier-critic-disturbance approximate optimal
  control of nonzero-sum differential games for modular robot manipulators}.
\newblock {\it \bibinfo{journal}{Neural Networks}\/},  {\it
  \bibinfo{volume}{181}\/}, \bibinfo{pages}{106880}.
\bibitem[{D’ambrosio et~al.(2021)D’ambrosio, Schiassi, Curti \&
  Furfaro}]{d2021pontryagin}
\bibinfo{author}{D’ambrosio, A.}, \bibinfo{author}{Schiassi, E.},
  \bibinfo{author}{Curti, F.}, \& \bibinfo{author}{Furfaro, R.}
  (\bibinfo{year}{2021}).
\newblock \bibinfo{title}{Pontryagin neural networks with functional
  interpolation for optimal intercept problems}.
\newblock {\it \bibinfo{journal}{Mathematics}\/},  {\it \bibinfo{volume}{9}\/},
  \bibinfo{pages}{996}.
\bibitem[{Engwerda(2005)}]{engwerda2005lq}
\bibinfo{author}{Engwerda, J.} (\bibinfo{year}{2005}).
\newblock {\it \bibinfo{title}{LQ dynamic optimization and differential
  games}\/}.
\newblock \bibinfo{publisher}{John Wiley \& Sons}.
\bibitem[{Furfaro et~al.(2022)Furfaro, D'Ambrosio, Schiassi \&
  Scorsoglio}]{furfaro2022physics}
\bibinfo{author}{Furfaro, R.}, \bibinfo{author}{D'Ambrosio, A.},
  \bibinfo{author}{Schiassi, E.}, \& \bibinfo{author}{Scorsoglio, A.}
  (\bibinfo{year}{2022}).
\newblock \bibinfo{title}{Physics-informed neural networks for closed-loop
  guidance and control in aerospace systems}.
\newblock In {\it \bibinfo{booktitle}{AIAA SCITECH 2022 Forum}\/} (p.
  \bibinfo{pages}{0361}).
\bibitem[{Hager(1990)}]{hager1990multiplier}
\bibinfo{author}{Hager, W.~W.} (\bibinfo{year}{1990}).
\newblock \bibinfo{title}{Multiplier methods for nonlinear optimal control}.
\newblock {\it \bibinfo{journal}{SIAM Journal on Numerical Analysis}\/},  {\it
  \bibinfo{volume}{27}\/}, \bibinfo{pages}{1061--1080}.
\bibitem[{Hannemann-Tam{\'a}s \& Marquardt(2012)}]{hannemann2012verify}
\bibinfo{author}{Hannemann-Tam{\'a}s, R.}, \& \bibinfo{author}{Marquardt, W.}
  (\bibinfo{year}{2012}).
\newblock \bibinfo{title}{How to verify optimal controls computed by direct
  shooting methods?--a tutorial}.
\newblock {\it \bibinfo{journal}{Journal of process control}\/},  {\it
  \bibinfo{volume}{22}\/}, \bibinfo{pages}{494--507}.
\bibitem[{He et~al.(2025)He, Li \& Nie}]{he2025Synchronization}
\bibinfo{author}{He, Z.}, \bibinfo{author}{Li, C.}, \& \bibinfo{author}{Nie,
  L.} (\bibinfo{year}{2025}).
\newblock \bibinfo{title}{Synchronization of complex dynamical networks with
  saturated delayed impulsive control}.
\newblock {\it \bibinfo{journal}{ISA Transactions}\/},  {\it
  \bibinfo{volume}{157}\/}, \bibinfo{pages}{153--163}.
\bibitem[{Hou et~al.(2023)Hou, Li \& Ying}]{hou2023enhancing}
\bibinfo{author}{Hou, J.}, \bibinfo{author}{Li, Y.}, \& \bibinfo{author}{Ying,
  S.} (\bibinfo{year}{2023}).
\newblock \bibinfo{title}{Enhancing pinns for solving pdes via adaptive
  collocation point movement and adaptive loss weighting}.
\newblock {\it \bibinfo{journal}{Nonlinear Dynamics}\/},  {\it
  \bibinfo{volume}{111}\/}, \bibinfo{pages}{15233--15261}.
\bibitem[{Hwang et~al.(2022)Hwang, Lee, Shin \& Hwang}]{hwang2022solving}
\bibinfo{author}{Hwang, R.}, \bibinfo{author}{Lee, J.~Y.},
  \bibinfo{author}{Shin, J.~Y.}, \& \bibinfo{author}{Hwang, H.~J.}
  (\bibinfo{year}{2022}).
\newblock \bibinfo{title}{Solving pde-constrained control problems using
  operator learning}.
\newblock In {\it \bibinfo{booktitle}{Proceedings of the AAAI Conference on
  Artificial Intelligence}\/} (pp. \bibinfo{pages}{4504--4512}).
\newblock volume~\bibinfo{volume}{36}.
\bibitem[{Ishihara \& Morimoto(2018)}]{ishihara2018optimal}
\bibinfo{author}{Ishihara, K.}, \& \bibinfo{author}{Morimoto, J.}
  (\bibinfo{year}{2018}).
\newblock \bibinfo{title}{An optimal control strategy for hybrid actuator
  systems: Application to an artificial muscle with electric motor assist}.
\newblock {\it \bibinfo{journal}{Neural Networks}\/},  {\it
  \bibinfo{volume}{99}\/}, \bibinfo{pages}{92--100}.
\bibitem[{Kendall et~al.(2018)Kendall, Gal \& Cipolla}]{kendall2018multi}
\bibinfo{author}{Kendall, A.}, \bibinfo{author}{Gal, Y.}, \&
  \bibinfo{author}{Cipolla, R.} (\bibinfo{year}{2018}).
\newblock \bibinfo{title}{Multi-task learning using uncertainty to weigh losses
  for scene geometry and semantics}.
\newblock In {\it \bibinfo{booktitle}{Proceedings of the IEEE conference on
  computer vision and pattern recognition}\/} (pp.
  \bibinfo{pages}{7482--7491}).
\bibitem[{Kirk(2004)}]{kirk2004optimal}
\bibinfo{author}{Kirk, D.~E.} (\bibinfo{year}{2004}).
\newblock {\it \bibinfo{title}{Optimal control theory: an introduction}\/}.
\newblock \bibinfo{publisher}{Courier Corporation}.
\bibitem[{Koumir et~al.(2016)Koumir, El~Bakri \& Boumhidi}]{koumir2016optimal}
\bibinfo{author}{Koumir, M.}, \bibinfo{author}{El~Bakri, A.}, \&
  \bibinfo{author}{Boumhidi, I.} (\bibinfo{year}{2016}).
\newblock \bibinfo{title}{Optimal control for a variable speed wind turbine
  based on extreme learning machine and adaptive particle swarm optimization}.
\newblock In {\it \bibinfo{booktitle}{2016 5th International Conference on
  Systems and Control (ICSC)}\/} (pp. \bibinfo{pages}{151--156}).
\newblock \bibinfo{organization}{IEEE}.
\bibitem[{Lasdon et~al.(2003)Lasdon, Mitter \& Waren}]{lasdon2003conjugate}
\bibinfo{author}{Lasdon, L.}, \bibinfo{author}{Mitter, S.}, \&
  \bibinfo{author}{Waren, A.} (\bibinfo{year}{2003}).
\newblock \bibinfo{title}{The conjugate gradient method for optimal control
  problems}.
\newblock {\it \bibinfo{journal}{IEEE Transactions on Automatic Control}\/},
  {\it \bibinfo{volume}{12}\/}, \bibinfo{pages}{132--138}.
\bibitem[{Li et~al.(2024{\natexlab{a}})Li, Tan, Su et~al.}]{Li2024Verification}
\bibinfo{author}{Li, R.}, \bibinfo{author}{Tan, Y.}, \bibinfo{author}{Su, X.}
  et~al. (\bibinfo{year}{2024}{\natexlab{a}}).
\newblock \bibinfo{title}{A verification theorem for feedback nash equilibrium
  in multiple-player nonzero-sum impulse game}.
\newblock {\it \bibinfo{journal}{IEEE/CAA Journal of Automatica Sinica}\/}, .
\newblock \bibinfo{note}{In press}.
\bibitem[{Li et~al.(2018)Li, Jusup, Wang, Li, Shi, Podobnik, Stanley, Havlin \&
  Boccaletti}]{li2018punishment}
\bibinfo{author}{Li, X.}, \bibinfo{author}{Jusup, M.}, \bibinfo{author}{Wang,
  Z.}, \bibinfo{author}{Li, H.}, \bibinfo{author}{Shi, L.},
  \bibinfo{author}{Podobnik, B.}, \bibinfo{author}{Stanley, H.~E.},
  \bibinfo{author}{Havlin, S.}, \& \bibinfo{author}{Boccaletti, S.}
  (\bibinfo{year}{2018}).
\newblock \bibinfo{title}{Punishment diminishes the benefits of network
  reciprocity in social dilemma experiments}.
\newblock {\it \bibinfo{journal}{Proceedings of the National Academy of
  Sciences}\/},  {\it \bibinfo{volume}{115}\/}, \bibinfo{pages}{30--35}.
\bibitem[{Li et~al.(2024{\natexlab{b}})Li, Sun, Marques, Wang \&
  You}]{li2024pontryagin}
\bibinfo{author}{Li, Z.}, \bibinfo{author}{Sun, J.}, \bibinfo{author}{Marques,
  A.~G.}, \bibinfo{author}{Wang, G.}, \& \bibinfo{author}{You, K.}
  (\bibinfo{year}{2024}{\natexlab{b}}).
\newblock \bibinfo{title}{Pontryagin’s minimum principle-guided rl for
  minimum-time exploration of spatiotemporal fields}.
\newblock {\it \bibinfo{journal}{IEEE Transactions on Neural Networks and
  Learning Systems}\/}, .
\bibitem[{Liu et~al.(2023{\natexlab{a}})Liu, Liao, Dong \&
  Mansoori}]{liu2023neurodynamic}
\bibinfo{author}{Liu, J.}, \bibinfo{author}{Liao, X.}, \bibinfo{author}{Dong,
  J.-s.}, \& \bibinfo{author}{Mansoori, A.}
  (\bibinfo{year}{2023}{\natexlab{a}}).
\newblock \bibinfo{title}{A neurodynamic approach for nonsmooth optimal power
  consumption of intelligent and connected vehicles}.
\newblock {\it \bibinfo{journal}{Neural Networks}\/},  {\it
  \bibinfo{volume}{161}\/}, \bibinfo{pages}{693--707}.
\bibitem[{Liu et~al.(2023{\natexlab{b}})Liu, Ding, Zhang \& Zhou}]{liu2023pinn}
\bibinfo{author}{Liu, T.}, \bibinfo{author}{Ding, S.}, \bibinfo{author}{Zhang,
  J.}, \& \bibinfo{author}{Zhou, L.} (\bibinfo{year}{2023}{\natexlab{b}}).
\newblock \bibinfo{title}{Pinn-based viscosity solution of hjb equation}.
\newblock {\it \bibinfo{journal}{arXiv preprint arXiv:2309.09953}\/}, .
\bibitem[{Liu et~al.(2024)Liu, Gu, Yu \& Qin}]{liu2024diminishing}
\bibinfo{author}{Liu, Y.}, \bibinfo{author}{Gu, H.}, \bibinfo{author}{Yu, X.},
  \& \bibinfo{author}{Qin, P.} (\bibinfo{year}{2024}).
\newblock \bibinfo{title}{Diminishing spectral bias in physics-informed neural
  networks using spatially-adaptive fourier feature encoding}.
\newblock {\it \bibinfo{journal}{Neural Networks}\/},  (p.
  \bibinfo{pages}{106886}).
\bibitem[{Lv \& Xiong(2022)}]{lv2022nonzero}
\bibinfo{author}{Lv, S.}, \& \bibinfo{author}{Xiong, J.}
  (\bibinfo{year}{2022}).
\newblock \bibinfo{title}{Nonzero-sum impulse games with regime switching}.
\newblock {\it \bibinfo{journal}{Automatica}\/},  {\it
  \bibinfo{volume}{145}\/}, \bibinfo{pages}{110439}.
\bibitem[{Mowlavi \& Nabi(2023)}]{mowlavi2023optimal}
\bibinfo{author}{Mowlavi, S.}, \& \bibinfo{author}{Nabi, S.}
  (\bibinfo{year}{2023}).
\newblock \bibinfo{title}{Optimal control of pdes using physics-informed neural
  networks}.
\newblock {\it \bibinfo{journal}{Journal of Computational Physics}\/},  {\it
  \bibinfo{volume}{473}\/}, \bibinfo{pages}{111731}.
\bibitem[{Mukherjee \& Liu(2023)}]{mukherjee2023bridging}
\bibinfo{author}{Mukherjee, A.}, \& \bibinfo{author}{Liu, J.}
  (\bibinfo{year}{2023}).
\newblock \bibinfo{title}{Bridging physics-informed neural networks with
  reinforcement learning: Hamilton-jacobi-bellman proximal policy optimization
  (hjbppo)}.
\newblock {\it \bibinfo{journal}{arXiv preprint arXiv:2302.00237}\/}, .
\bibitem[{Nikooeinejad et~al.(2016)Nikooeinejad, Delavarkhalafi \&
  Heydari}]{nikooeinejad2016numerical}
\bibinfo{author}{Nikooeinejad, Z.}, \bibinfo{author}{Delavarkhalafi, A.}, \&
  \bibinfo{author}{Heydari, M.} (\bibinfo{year}{2016}).
\newblock \bibinfo{title}{A numerical solution of open-loop nash equilibrium in
  nonlinear differential games based on chebyshev pseudospectral method}.
\newblock {\it \bibinfo{journal}{Journal of Computational and Applied
  Mathematics}\/},  {\it \bibinfo{volume}{300}\/}, \bibinfo{pages}{369--384}.
\bibitem[{Paszke et~al.(2019)Paszke, Gross, Massa, Lerer, Bradbury, Chanan,
  Killeen, Lin, Gimelshein, Antiga et~al.}]{paszke2019pytorch}
\bibinfo{author}{Paszke, A.}, \bibinfo{author}{Gross, S.},
  \bibinfo{author}{Massa, F.}, \bibinfo{author}{Lerer, A.},
  \bibinfo{author}{Bradbury, J.}, \bibinfo{author}{Chanan, G.},
  \bibinfo{author}{Killeen, T.}, \bibinfo{author}{Lin, Z.},
  \bibinfo{author}{Gimelshein, N.}, \bibinfo{author}{Antiga, L.} et~al.
  (\bibinfo{year}{2019}).
\newblock \bibinfo{title}{Pytorch: An imperative style, high-performance deep
  learning library}.
\newblock {\it \bibinfo{journal}{Advances in neural information processing
  systems}\/},  {\it \bibinfo{volume}{32}\/}.
\bibitem[{Qiu et~al.(2022)Qiu, Li \& Qiu}]{Qiu2022Fuzzy}
\bibinfo{author}{Qiu, R.}, \bibinfo{author}{Li, R.}, \& \bibinfo{author}{Qiu,
  J.} (\bibinfo{year}{2022}).
\newblock \bibinfo{title}{A novel step-function method for stability analysis
  of t-s fuzzy impulsive systems}.
\newblock {\it \bibinfo{journal}{IEEE Transactions on Fuzzy Systems}\/},  {\it
  \bibinfo{volume}{30}\/}, \bibinfo{pages}{4399--4408}.
\bibitem[{Raissi et~al.(2019)Raissi, Perdikaris \&
  Karniadakis}]{raissi2019physics}
\bibinfo{author}{Raissi, M.}, \bibinfo{author}{Perdikaris, P.}, \&
  \bibinfo{author}{Karniadakis, G.~E.} (\bibinfo{year}{2019}).
\newblock \bibinfo{title}{Physics-informed neural networks: A deep learning
  framework for solving forward and inverse problems involving nonlinear
  partial differential equations}.
\newblock {\it \bibinfo{journal}{Journal of Computational physics}\/},  {\it
  \bibinfo{volume}{378}\/}, \bibinfo{pages}{686--707}.
\bibitem[{Robinson et~al.(2022)Robinson, Pawar, Rasheed \&
  San}]{robinson2022physics}
\bibinfo{author}{Robinson, H.}, \bibinfo{author}{Pawar, S.},
  \bibinfo{author}{Rasheed, A.}, \& \bibinfo{author}{San, O.}
  (\bibinfo{year}{2022}).
\newblock \bibinfo{title}{Physics guided neural networks for modelling of
  non-linear dynamics}.
\newblock {\it \bibinfo{journal}{Neural Networks}\/},  {\it
  \bibinfo{volume}{154}\/}, \bibinfo{pages}{333--345}.
\bibitem[{Sadana et~al.(2022)Sadana, Reddy \& Zaccour}]{sadana2022feedback}
\bibinfo{author}{Sadana, U.}, \bibinfo{author}{Reddy, P.~V.}, \&
  \bibinfo{author}{Zaccour, G.} (\bibinfo{year}{2022}).
\newblock \bibinfo{title}{Feedback nash equilibria in differential games with
  impulse control}.
\newblock {\it \bibinfo{journal}{IEEE Transactions on Automatic Control}\/},
  {\it \bibinfo{volume}{68}\/}, \bibinfo{pages}{4523--4538}.
\bibitem[{Sager(2009)}]{sager2009reformulations}
\bibinfo{author}{Sager, S.} (\bibinfo{year}{2009}).
\newblock \bibinfo{title}{Reformulations and algorithms for the optimization of
  switching decisions in nonlinear optimal control}.
\newblock {\it \bibinfo{journal}{Journal of Process Control}\/},  {\it
  \bibinfo{volume}{19}\/}, \bibinfo{pages}{1238--1247}.
\bibitem[{Schiassi et~al.(2022)Schiassi, D’Ambrosio, Drozd, Curti \&
  Furfaro}]{schiassi2022physics}
\bibinfo{author}{Schiassi, E.}, \bibinfo{author}{D’Ambrosio, A.},
  \bibinfo{author}{Drozd, K.}, \bibinfo{author}{Curti, F.}, \&
  \bibinfo{author}{Furfaro, R.} (\bibinfo{year}{2022}).
\newblock \bibinfo{title}{Physics-informed neural networks for optimal planar
  orbit transfers}.
\newblock {\it \bibinfo{journal}{Journal of Spacecraft and Rockets}\/},  {\it
  \bibinfo{volume}{59}\/}, \bibinfo{pages}{834--849}.
\bibitem[{Schiassi et~al.(2020)Schiassi, D’Ambrosio, Johnston, De~Florio,
  Drozd, Furfaro, Curti, Mortari et~al.}]{schiassi2020physics}
\bibinfo{author}{Schiassi, E.}, \bibinfo{author}{D’Ambrosio, A.},
  \bibinfo{author}{Johnston, H.}, \bibinfo{author}{De~Florio, M.},
  \bibinfo{author}{Drozd, K.}, \bibinfo{author}{Furfaro, R.},
  \bibinfo{author}{Curti, F.}, \bibinfo{author}{Mortari, D.} et~al.
  (\bibinfo{year}{2020}).
\newblock \bibinfo{title}{Physics-informed extreme theory of functional
  connections applied to optimal orbit transfer}.
\newblock In {\it \bibinfo{booktitle}{Proceedings of the AAS/AIAA Astrodynamics
  Specialist Conference, Lake Tahoe, CA, USA}\/} (pp. \bibinfo{pages}{9--13}).
\bibitem[{Schiassi et~al.(2021)Schiassi, D’Ambrosio, Scorsoglio, Furfaro \&
  Curti}]{schiassi2021class}
\bibinfo{author}{Schiassi, E.}, \bibinfo{author}{D’Ambrosio, A.},
  \bibinfo{author}{Scorsoglio, A.}, \bibinfo{author}{Furfaro, R.}, \&
  \bibinfo{author}{Curti, F.} (\bibinfo{year}{2021}).
\newblock \bibinfo{title}{Class of optimal space guidance problems solved via
  indirect methods and physics-informed neural networks}.
\newblock In {\it \bibinfo{booktitle}{Proceedings of 31st AAS/AIAA Space Flight
  Mechanics Meeting}\/}.
\bibitem[{Shilova et~al.(2024)Shilova, Delliaux, Preux \&
  Raffin}]{shilova2024learning}
\bibinfo{author}{Shilova, A.}, \bibinfo{author}{Delliaux, T.},
  \bibinfo{author}{Preux, P.}, \& \bibinfo{author}{Raffin, B.}
  (\bibinfo{year}{2024}).
\newblock {\it \bibinfo{title}{Learning HJB Viscosity Solutions with PINNs for
  Continuous-Time Reinforcement Learning}\/}.
\newblock Ph.D. thesis Inria Lille-Nord Europe, CRIStAL-Centre de Recherche en
  Informatique, Signal et Automatique de Lille-UMR 9189.
\bibitem[{Tauchnitz(2015)}]{tauchnitz2015pontryagin}
\bibinfo{author}{Tauchnitz, N.} (\bibinfo{year}{2015}).
\newblock \bibinfo{title}{The pontryagin maximum principle for nonlinear
  optimal control problems with infinite horizon}.
\newblock {\it \bibinfo{journal}{Journal of Optimization Theory and
  Applications}\/},  {\it \bibinfo{volume}{167}\/}, \bibinfo{pages}{27--48}.
\bibitem[{Wang et~al.(2024{\natexlab{a}})Wang, Pang, Xu, Huang \&
  Kurths}]{Wang2024State}
\bibinfo{author}{Wang, X.}, \bibinfo{author}{Pang, N.}, \bibinfo{author}{Xu,
  Y.}, \bibinfo{author}{Huang, T.}, \& \bibinfo{author}{Kurths, J.}
  (\bibinfo{year}{2024}{\natexlab{a}}).
\newblock \bibinfo{title}{On state constrained containment control for
  nonlinear multiagent systems using event-triggered input}.
\newblock {\it \bibinfo{journal}{IEEE Transactions on Systems, Man, and
  Cybernetics: Systems}\/},  {\it \bibinfo{volume}{54}\/},
  \bibinfo{pages}{2530--2538}.
\bibitem[{Wang et~al.(2024{\natexlab{b}})Wang, Xu, Huang \&
  Kurths}]{Wang2024EventTriggered}
\bibinfo{author}{Wang, X.}, \bibinfo{author}{Xu, R.}, \bibinfo{author}{Huang,
  T.}, \& \bibinfo{author}{Kurths, J.} (\bibinfo{year}{2024}{\natexlab{b}}).
\newblock \bibinfo{title}{Event-triggered adaptive containment control for
  heterogeneous stochastic nonlinear multiagent systems}.
\newblock {\it \bibinfo{journal}{IEEE Transactions on Neural Networks and
  Learning Systems}\/},  {\it \bibinfo{volume}{35}\/},
  \bibinfo{pages}{8524--8534}.
\bibitem[{Wang et~al.(2024{\natexlab{c}})Wang, Li, Wu \&
  Deng}]{wang2024Stability}
\bibinfo{author}{Wang, Y.}, \bibinfo{author}{Li, C.}, \bibinfo{author}{Wu, H.},
  \& \bibinfo{author}{Deng, H.} (\bibinfo{year}{2024}{\natexlab{c}}).
\newblock \bibinfo{title}{Stability of nonlinear delayed impulsive control
  systems via step-function method}.
\newblock {\it \bibinfo{journal}{Chaos, Solitons and Fractals: the
  interdisciplinary journal of Nonlinear Science, and Nonequilibrium and
  Complex Phenomena}\/},  {\it \bibinfo{volume}{189}\/}.
\bibitem[{Wang et~al.(2025)Wang, Song \& Liu}]{wang2025synchronisation}
\bibinfo{author}{Wang, Y.}, \bibinfo{author}{Song, Q.}, \&
  \bibinfo{author}{Liu, Y.} (\bibinfo{year}{2025}).
\newblock \bibinfo{title}{Synchronisation of quaternion-valued neural networks
  with neutral delay and discrete delay via aperiodic intermittent control}.
\newblock {\it \bibinfo{journal}{International Journal of Systems Science}\/},
  {\it \bibinfo{volume}{56}\/}, \bibinfo{pages}{1395--1412}.
\bibitem[{Wang et~al.(2018)Wang, Jusup, Shi, Lee, Iwasa \&
  Boccaletti}]{wang2018exploiting}
\bibinfo{author}{Wang, Z.}, \bibinfo{author}{Jusup, M.}, \bibinfo{author}{Shi,
  L.}, \bibinfo{author}{Lee, J.-H.}, \bibinfo{author}{Iwasa, Y.}, \&
  \bibinfo{author}{Boccaletti, S.} (\bibinfo{year}{2018}).
\newblock \bibinfo{title}{Exploiting a cognitive bias promotes cooperation in
  social dilemma experiments}.
\newblock {\it \bibinfo{journal}{Nature communications}\/},  {\it
  \bibinfo{volume}{9}\/}, \bibinfo{pages}{2954}.
\bibitem[{Wang et~al.(2017)Wang, Jusup, Wang, Shi, Iwasa, Moreno \&
  Kurths}]{wang2017onymity}
\bibinfo{author}{Wang, Z.}, \bibinfo{author}{Jusup, M.}, \bibinfo{author}{Wang,
  R.-W.}, \bibinfo{author}{Shi, L.}, \bibinfo{author}{Iwasa, Y.},
  \bibinfo{author}{Moreno, Y.}, \& \bibinfo{author}{Kurths, J.}
  (\bibinfo{year}{2017}).
\newblock \bibinfo{title}{Onymity promotes cooperation in social dilemma
  experiments}.
\newblock {\it \bibinfo{journal}{Science advances}\/},  {\it
  \bibinfo{volume}{3}\/}, \bibinfo{pages}{e1601444}.
\bibitem[{Wang et~al.(2015{\natexlab{a}})Wang, Kokubo, Jusup \&
  Tanimoto}]{wang2015universal}
\bibinfo{author}{Wang, Z.}, \bibinfo{author}{Kokubo, S.},
  \bibinfo{author}{Jusup, M.}, \& \bibinfo{author}{Tanimoto, J.}
  (\bibinfo{year}{2015}{\natexlab{a}}).
\newblock \bibinfo{title}{Universal scaling for the dilemma strength in
  evolutionary games}.
\newblock {\it \bibinfo{journal}{Physics of life reviews}\/},  {\it
  \bibinfo{volume}{14}\/}, \bibinfo{pages}{1--30}.
\bibitem[{Wang et~al.(2022)Wang, Mu, Hu, Chu \& Li}]{wang2022modelling}
\bibinfo{author}{Wang, Z.}, \bibinfo{author}{Mu, C.}, \bibinfo{author}{Hu, S.},
  \bibinfo{author}{Chu, C.}, \& \bibinfo{author}{Li, X.}
  (\bibinfo{year}{2022}).
\newblock \bibinfo{title}{Modelling the dynamics of regret minimization in
  large agent populations: a master equation approach.}
\newblock In {\it \bibinfo{booktitle}{IJCAI}\/} (pp.
  \bibinfo{pages}{534--540}).
\bibitem[{Wang et~al.(2015{\natexlab{b}})Wang, Wang, Szolnoki \&
  Perc}]{wang2015evolutionary}
\bibinfo{author}{Wang, Z.}, \bibinfo{author}{Wang, L.},
  \bibinfo{author}{Szolnoki, A.}, \& \bibinfo{author}{Perc, M.}
  (\bibinfo{year}{2015}{\natexlab{b}}).
\newblock \bibinfo{title}{Evolutionary games on multilayer networks: a
  colloquium}.
\newblock {\it \bibinfo{journal}{The European physical journal B}\/},  {\it
  \bibinfo{volume}{88}\/}, \bibinfo{pages}{1--15}.
\bibitem[{Wu \& Lisser(2023)}]{wu2023enhancing}
\bibinfo{author}{Wu, D.}, \& \bibinfo{author}{Lisser, A.}
  (\bibinfo{year}{2023}).
\newblock \bibinfo{title}{Enhancing neurodynamic approach with physics-informed
  neural networks for solving non-smooth convex optimization problems}.
\newblock {\it \bibinfo{journal}{Neural Networks}\/},  {\it
  \bibinfo{volume}{168}\/}, \bibinfo{pages}{419--430}.
\bibitem[{Xiang et~al.(2022)Xiang, Peng, Liu \& Yao}]{xiang2022self}
\bibinfo{author}{Xiang, Z.}, \bibinfo{author}{Peng, W.}, \bibinfo{author}{Liu,
  X.}, \& \bibinfo{author}{Yao, W.} (\bibinfo{year}{2022}).
\newblock \bibinfo{title}{Self-adaptive loss balanced physics-informed neural
  networks}.
\newblock {\it \bibinfo{journal}{Neurocomputing}\/},  {\it
  \bibinfo{volume}{496}\/}, \bibinfo{pages}{11--34}.
\bibitem[{Yang et~al.(2024)Yang, Yang, Deng \& He}]{yang2024moving}
\bibinfo{author}{Yang, Y.}, \bibinfo{author}{Yang, Q.}, \bibinfo{author}{Deng,
  Y.}, \& \bibinfo{author}{He, Q.} (\bibinfo{year}{2024}).
\newblock \bibinfo{title}{Moving sampling physics-informed neural networks
  induced by moving mesh pde}.
\newblock {\it \bibinfo{journal}{Neural Networks}\/},  {\it
  \bibinfo{volume}{180}\/}, \bibinfo{pages}{106706}.
\bibitem[{Zhang et~al.(2024)Zhang, Ghimire, Zhang, Xu \& Ren}]{zhang2024value}
\bibinfo{author}{Zhang, L.}, \bibinfo{author}{Ghimire, M.},
  \bibinfo{author}{Zhang, W.}, \bibinfo{author}{Xu, Z.}, \&
  \bibinfo{author}{Ren, Y.} (\bibinfo{year}{2024}).
\newblock \bibinfo{title}{Value approximation for two-player general-sum
  differential games with state constraints}.
\newblock {\it \bibinfo{journal}{IEEE Transactions on Robotics}\/}, .

\end{thebibliography}

\end{document}